\documentclass[12pt, dvips, twoside]{amsart}
\usepackage{amsmath,amssymb,amscd,amsxtra,latexsym,bbm,graphicx,epsfig,epic,
eepic,oldgerm,bm,psfrag,amsthm}
\usepackage[mathscr]{euscript}
\usepackage[active]{srcltx}
\usepackage{a4wide}

\input xy
\xyoption{all}
\CompileMatrices
    
\theoremstyle{plain}
\newtheorem{theorem}{Theorem}[section]
\newtheorem*{theorem*}{Theorem}
\newtheorem{lemma}[theorem]{Lemma}
\newtheorem{proposition}[theorem]{Proposition}

\newtheorem*{remark*}{Remark}
\newtheorem*{remarks*}{Remarks}

\newtheorem{example}[theorem]{Example}
\newtheorem*{example*}{Example}
\newtheorem*{examples*}{Examples}

\newtheorem*{definition*}{Definition}

\newcommand*{\lo}[1]{{\raisebox{-.4ex}{\mathsurround=0pt\makebox{$\scriptstyle\mskip-.8mu#1$}}}}

\newcommand{\proofend}{\hspace*{\fill} $\Box$\\}
\newcommand{\diam}{\hspace*{\fill} $\Diamond$\\}

\newcommand{\ign}[1]{}

\def\1{\:\!}
\def\2{\;\!}

\def\s{\smallskip}
\def\m{\medskip}

\def\Diffc0{\operatorname{Diff^c_0}}
\def\Symp{\operatorname{Symp}}

\def\Sympc0{\operatorname{Symp^c_0}}

\def\Ham{\operatorname{Ham}}

\def\GL{\operatorname{GL}}

\def\Emb{\operatorname{Emb}}
\def\Sp{\operatorname{Sp}}

\def\rk{\operatorname{rk}}

\def\eps{\varepsilon}
\def\gve{\varepsilon}
\def\gf{\varphi}

\def\gl{\lambda}
\def\go{\omega}
\def\gs{\sigma}

\def\cc{{\mathcal C}}
\def\cd{{\mathscr D}}
\def\ce{{\mathcal E}}

\def\ch{{\mathcal H}}

\def\cj{{\mathcal J}}
\def\ck{{\mathcal K}}
\def\cl{{\mathcal L}}

\def\cv{{\mathcal V}}

\def\ua{\underline{\bm {a\?}}\.\.}

\def\ba{{\bm {a}}}
\def\bd{{\bm {d}}}
\def\bc{{\bm {c}}}
\def\be{{\bm {e}}}

\def\bn{{\bm {n}}}
\def\bp{{\bm {p}}}
\def\bq{{\bm {q}}}
\def\bs{{\bm {s}}}
\def\bv{{\bm {v}}}
\def\bu{{\bm {u}}}
\def\bw{{\bm {w}}}

\def\T{\widehat{\mbox{$\mskip-.2mu T\mskip1mu$}}}
\def\t{\widehat{\mbox{$\mskip.5mu\scriptstyle\mskip-.2mu T\mskip.5mu$}}}

\def\CC{\mathbbm{C}}

\def\NN{\mathbbm{N}}

\def\QQ{\mathbbm{Q}}
\def\RR{\mathbbm{R}}

\def\ZZ{\mathbbm{Z}}

\def\CP{\operatorname{\mathbbm{C}P}}

\def\pp{\partial}

\def\ni{\noindent}
\def\b{\bigskip}
\def\m{\medskip}

\def\.{\mskip1mu}
\def\?{\mskip-1mu}

\newcommand{\Bcirc}{\overset%
{\raisebox{-.3ex}[0ex][-.3ex]{\mbox{$\scriptscriptstyle\circ$}}\mskip-5mu}B{}}

\newcommand{\gdot}{\dot {\gamma\mskip2.5mu}\mskip-2.5mu}

\def\id{\operatorname{\mbox{\tt id}}}

\def\proof{\noindent {\it Proof. \;}}
\newcommand{\proofof}[1]{\ni {\it Proof of #1. }}

\begin{document}

\title{Lagrangian product tori in tame symplectic manifolds}

\author{Yuri Chekanov} \thanks{YuC partially supported by RFBR grant NSh-5138.2014.1}
\address{(Yu.~Chekanov)
Moscow Center for Continuous Mathematical Education,
B.~Vlasievsky per.~11, Moscow 121002, Russia}
\email{chekanov@mccme.ru}
\author{Felix Schlenk} \thanks{FS partially supported by SNF grant 200020-144432/1.}
\address{(F.~Schlenk)
Institut de Math\'ematiques,
Universit\'e de Neuch\^atel,
Rue \'Emile Argand~11,
CP~158,
2009 Neuch\^atel,
Switzerland}
\email{schlenk@unine.ch}

\date{\today}

\begin{abstract}
In \cite{Ch-tori}, product Lagrangian tori in standard symplectic
space~$\RR^{2n}$ were classified up to symplectomorphism.
We extend this classification to tame symplectically aspherical symplectic manifolds.
We show by examples that the asphericity assumption cannot be omitted.
\end{abstract}

\maketitle

\markboth{{\rm Lagrangian product tori in tame symplectic manifolds}}{{}}

%\begin{center}
%{\it
%Dedicated to Edi Zehnder 
%at the occasion of his 75th birthday
%}
%\end{center}

\section{Introduction and main results}

\ni
Let $T(a)$ denote the boundary of the disc of area $a>0$ in~$\RR^2$
centred at the origin.
Let $\ba \.=\.(a_1,\ldots,a_n)$ be a vector  with positive components.
We call the $n$-torus
$$
T(\ba) \,=\, T(a_1) \times \dots \times T(a_n) \subset \RR^{2n}
$$
a {\it product torus}.
Product tori are Lagrangian with respect to the standard symplectic form
$\.\.\go_n=\sum_{j=1}^n dx_j \wedge dy_j$, that is,
the restriction of~$\go_n$ to each product torus vanishes.

Let $(M,\go)$ be a symplectic manifold.
We assume throughout the paper that $M$ is connected.
Denote by $B^{2n}(b)$ the closed ball of radius $\sqrt{b/\pi}$ in
$\RR^{2n}$ centred at the origin.
The torus $T(\ba)$ lies on the boundary of the ball~$B^{2n}(|\ba|)$.
By a {\it symplectic chart\/} we understand a symplectic embedding
$\gf \colon B^{2n}(b)\to (M,\go)$.
%F
%(This means that there is an neighbourhood~$U$ of $B^{2n}(b)$ in~$\RR^{2n}$
%and a symplectic embedding $\gf \colon U \to M$.)
%
Given  a symplectic chart $\gf \colon B^{2n}(b) \to (M,\go)$
and a torus $T(\ba) \subset B^{2n}(b)$,
we write  $T_\gf(\ba) = \gf \left( T(\ba) \right)$.
A Lagrangian torus in $(M,\go)$ is called a {\it product torus}\/
if it is of the form $T_\gf(\ba)$ for some symplectic chart~$\gf$.

We study the classification problem for product Lagrangian tori
with respect to the action of the group
$\Symp (M,\go)$ of {\it symplectomorphisms\/} of~$M$
(diffeomorphisms preserving the symplectic form~$\go$)
as well as the group $\Ham (M,\go)$
of {\it Hamiltonian symplectomorphisms}.
Hamiltonian symplectomorphisms are defined as follows.
Let $\{H_t\}$ be a family of smooth functions on $M$
smoothly depending on the parameter $t\in[0,1]$.
This family
defines a family of Hamiltonian vector fields  $\{X_{t}\}$
by $\go (X_{t}, \cdot \.) = dH_t (\.\cdot\.)$.
Assume that the time $t$ flow $\Psi_t$ of $\{X_{t}\}$
is a well-defined diffeomorphism for each $t\in[0,1]$.
Then each $\Psi_t$ is a symplectomorphism.
The family  $\{\Psi_t\}$ is then called a {\it Hamiltonian isotopy};
symplectomorphisms $\Psi_t$ arising in this way
form the subgroup $\Ham (M,\go) \subset \Symp (M,\go)$.

Given Lagrangian submanifolds $L, L'$
in a symplectic manifold~$(M,\go)$,
we write
$L \sim L'$ (resp.~$L \approx L'$)
if there is a symplectomorphism (resp.~a Hamiltonian symplectomorphism)
of $(M,\go)$ that maps $L$ to~$L'$.
In the particular case where $(M,\go)=(\RR^{2n},\go_n)$,
we say that~$L$ is Hamiltonian isotopic to~$L'$
in the ball~$B^{2n}(b)$
if there is a Hamiltonian isotopy $\{\Phi_s\}$, $s \in [0,1]$,
of~$\RR^{2n}$
such that $\Phi_0 = \id$, $\Phi_1(L) = L'$, and
$\Phi_s(L) \subset B^{2n}(b)$ for all $s \in [0,1]$.

Given a vector $\ba \.=\.(a_1,\ldots,a_n)$ with positive components, denote
\[
\ua = \min_{1\le i \le n} (a_i),
\,\,\,
m (\ba) = \# \left\{ i \mid a_i = \ua \right\},
\,\,\,
|\ba| = \textstyle\sum\limits_{i=1}^n  a_i,
\,\,\,
\|\ba\| = |\ba| + \ua.
\]
Let  $\Gamma (\ba)$ denote the subgroup of $\RR$ formed by all
integer combinations of  the numbers $a_1-\ua, \dots, a_n-\ua$.
We write $\ba\simeq  \ba'$ when the following holds:
$\ua=\ua'$, $m (\ba)=m (\ba')$, and $\Gamma (\ba)=\Gamma (\ba')$.
It was proved in \cite{Ch-tori} that for product tori
in $\RR^{2n}$ the conditions $T(\ba)\sim T(\ba')$, $T(\ba)\approx T(\ba')$,
$\ba\simeq \ba'$ are equivalent one to another.
The following theorem gives an upper bound on the size of the support of
a Hamiltonian isotopy between product tori when
such an isotopy exists.

\begin{theorem}  \label{t:t1}
{\rm (i)}
If $\ba$ and $\ba'$ are related by a permutation of the components,
then the tori $T(\ba)$ and  $T(\ba')$
are Hamiltonian isotopic in the ball $B^{2n}(|\ba|)$.

\s
{\rm (ii)}
If  $\ba\simeq\ba'$, then the tori $T(\ba)$ and  $T(\ba')$
are Hamiltonian isotopic in the ball $B^{2n}\bigl(\max(\|\ba\|,\|\ba'\|)\bigr)$.
\end{theorem}

Assertion~(i) of the theorem is, of course, rather obvious.
It seems likely that Theorem~\ref{t:t1} gives a sharp bound for the ball size.
However, we can only prove the sharpness  under the
additional assumption that $|\ba|\ne|\ba'|$:

\begin{theorem}  \label{t:t1+}
If $\,b<\max(\|\ba\|,\|\ba'\|)$ and $\,|\ba| \ne |\ba'|$,
then the tori $T(\ba)$ and $T(\ba')$
are not Hamiltonian isotopic in the ball~$B^{2n}(b)$.
\end{theorem}

%Along the proof of this result,
%we give obstructions to Hamiltonian isotopies between generalized Clifford tori
%in complex projective space~$\CP^n$, see Section~\ref{s:clifford}.

\m
It will sometimes be necessary  to assume that the geometry of the symplectic
manifold $\left( M,\go \right)$ is not too wild.
Following~\cite{Gr,Si,ALP}, we say that $\left( M,\go \right)$
is {\it tame}\/
if $M$ admits an almost complex structure $J$ and a complete
Riemannian metric $g$ satisfying the following conditions:
\begin{itemize}
\item[(T1)]
$J$ is uniformly tame, i.e.,
there are positive constants $C_1$ and $C_2$ such that
\[
\go \left( X, JX \right) \ge C_1 \left\| X \right\|_g^2
\quad \text{ and } \quad
\left| \.\.\go \left( X,Y \right) \right| \le C_2 \left\| X \right\|_g
\left\| Y \right\|_g
\]
for all tangent vectors $X$ and $Y $ on~$M$.

\s
\item[(T2)]
The sectional curvature of $\left(M,g\right)$ is bounded from above and
the injectivity radius of $\left(M,g\right)$ is bounded away from zero.
\end{itemize}

\ni
Some examples of tame symplectic manifolds are as follows:
\,(1) closed symplectic manifolds;
\,(2) cotangent bundles
over arbitrary  manifolds;
\,(3) twisted cotangent bundles over closed manifolds;
\,(4) symplectic manifolds such that
the complement of a compact subset is symplectomorphic to
the convex end of the symplectization of a closed contact manifold.
The class of tame symplectic manifolds is closed under taking
products and coverings.

\m
Recall that $(M,\go)$ is called {\it symplectically aspherical}\/ if 
$[\go] |_{\pi_2(M)} =0$ and $c_1 |_{\pi_2(M)}=0$.
Here, $c_1 = c_1(\go)$ is the first Chern class of~$TM$ with
respect to an (arbitrary) almost complex structure $J$ taming~$\go$ as in~(T1),
and the restriction to $\pi_2(M)$ is understood as the restriction to
the image of the natural map $\pi_2(M) \to H_2(M;\ZZ) \subset H_2(M;\RR)$.

Given  a symplectic chart $\gf \colon B^{2n}(b) \to (M,\go)$,
we write $b_\gf=b$.
The following theorem shows that the invariants
of product tori in $\RR^{2n}$ extend to
certain other symplectic manifolds:

\begin{theorem}  \label{t:tame}
Assume that $T_\gf(\ba)\sim T_{\gf'}(\ba')$,
where $T_\gf(\ba),  T_{\gf'}(\ba')\subset (M,\go)$.

\s
\begin{itemize}
\item[(i)]
If $(M,\go)$ is symplectically aspherical, then
$\Gamma (\ba)=\Gamma (\ba')$.

\s
\item[(ii)]
If $(M,\go)$ is tame, $\|\ba\| \le b_{\gf}$, and $\|\ba'\| \le b_{\gf'}$,
then $\ua = \ua'$ and $m(\ba) = m(\ba')$.
\end{itemize}
\end{theorem}

%\ni
%The assumption that $\Emb \left( B^{2n}(b),M,\go \right)$
%is connected cannot be omitted,
%as we shall show in Example~\ref{ex:camel}
%in Section~\ref{s:t2}.

\m
A symplectic manifold $(M,\go)$ is called a {\it Liouville manifold}\/
if it admits a vector field~$X$ such that $\cl_X \go  = \go$
(where $\cl_X$ is the Lie derivative with respect to~$X$).
If $X$ can be chosen in such a way that its time $t$ flow map is
well-defined for each $t \ge 0$, we call $(M,\go)$ {\it forward complete}.
Examples of tame forward complete Liouville manifolds are cotangent bundles and,
more generally, Stein manifolds, see~\cite{EG}.
Product tori in such manifolds can be  completely classified:

\begin{theorem}  \label{t:Liouville}
Let $T_\gf(\ba)$, $T_{\gf'}(\ba')$ be Lagrangian product tori
in a tame forward complete Liouville manifold~$(M,\go)$.
Then the conditions $\ba\simeq\ba'$, $T_\gf(\ba)\sim T_{\gf'}(\ba')$,
$T_\gf(\ba)\approx T_{\gf'}(\ba')$ are equivalent one to another.
\end{theorem}

The assumption $\| \ba \| \le b_\gf$, $\| \ba' \| \le b_{\gf'}$ in Theorem~\ref{t:tame}~(ii)
cannot be omitted, as the following simple example shows.
Let~$S^2$ be the round $2$-sphere, endowed with the Euclidean area form of
total area~$2$.
Let $p_N, p_S \in S^2$ be the north pole and the south pole.
Choose $\gve \in \;]0,\frac 12[$, and let $\gf, \gf' \colon B^2(2-\gve) \to S^2$
be Darboux charts such that $\gf(0) = p_N$, $\gf'(0) = p_S$, and such that
concentric circles are mapped to circles of latitude.
Then $T_{\gf}(\frac 12) = T_{\gf'}(\frac 32)$, but $\ua = \frac 12 \neq \frac 32 = \ua'$.
Note that $\| \ba' \| =3 >2-\gve = b_{\gf'}$.

\s
The assumption in Theorem~\ref{t:tame}~(i) that $(M,\go)$ is symplectically aspherical
cannot be omitted either, as the next theorem shows.
Recall that the cohomology class $[\go]$ of the symplectic form 
gives rise to the homomorphism
$\gs \colon \pi_2(M) \to \RR$, 
and  the first Chern class~$c_1 $ gives rise 
to the homomorphism $c_1 \colon \pi_2(M) \to \ZZ$.
Given $a>0$, define the homomorphism 
$$
\gs_a \colon \pi_2(M) \to \RR, \quad
S \mapsto \gs(S) - c_1(S) a .
$$
With $a>0$ 
and a symplectic manifold $(M,\go)$ we associate
the group 
$$
G_a\.=\.G_a(M,\go)\.:=\.\gs_a(\pi_2(M))\. \subset\. \RR.
$$
Note that $(M,\go)$ is symplectically aspherical if and only if~$G_a$
is trivial for all~$a>0$.
We call the symplectic manifold $(M,\go)$ {\it special}
if the rank of the group $\sigma\bigl(\pi_2(M)\bigr)\subset \RR$ is~$1$
and $c_1$ is not proportional to~$\sigma$.
We associate with each $S_0\in \pi_2(M)$ and each $a>0$
the  subgroup $G_a(S_0)=G_a(S_0,M,\go)$ of $G_a$ which is the image
under $\gs_a$ of the subgroup generated by~$S_0$.

\begin{theorem}  \label{t:shift}
Let $(M,\go)$ be a symplectic manifold;
if  $(M,\go)$ is special, we also fix an element $S_0\in \pi_2(M)$.
Let $\gf \colon B^{2n}(b) \to (M,\go)$ be a symplectic chart. 
For every real number~$c>0$ there exists $A>0$ such that
for all $a \in \:]0,A]$ the following holds.

If 
$d_1, \dots, d_k$ and $e_1, \dots, e_k$ for all  $j\in \{1, \dots, k\} $ satisfy 
the conditions
$d_j \ge c$,   $e_j \ge c$, 
$$
d_j-e_j \,\in\, 
\begin{cases} 
G_a (S_0)\,\,\,&\mbox{if $(M,\go)$ is special,}\\
G_a     \,\,\,&\mbox{otherwise},
\end{cases}
$$
and the tori $T_\gf (a, \dots, a, a+d_1, \dots, a+d_k )$,
$T_\gf (a, \dots, a, a+e_1, \dots, a+e_k )$ are contained  in $B_\gf(b)$,
then
\[
T_\gf (a, \dots, a, a+d_1, \dots, a+d_k ) 
\,\approx\, 
T_\gf (a, \dots, a, a+e_1, \dots, a+e_k ) .
\]
\end{theorem}

\begin{example*}
{\rm
Given $v>0$, we denote by $S^2(v)$ the $2$-sphere of area~$v$.
There exists a symplectic embedding 
$B^{4}(b) \to S^2(v_1) \times  S^2(v_2)$
whenever $b < \min ( v_1,  v_2 )$.
The homomorphism~$c_1$ on
$\pi_2(S^2(v_1) \times  S^2(v_2)) = \ZZ  \oplus \ZZ$ is given by
$(m_1,m_2) \mapsto 2 (m_1+m_2)$.
For $S_0=(1,-1)$ we have
$$
G_a\bigl(S_0, S^2(v_1) \times  S^2(v_2)\bigr) \,=\, 
 (v_1-v_2)  \.\.\ZZ .
$$
(Note that $S^2(v_1) \times  S^2(v_2)$ is special if and only if 
$v_1/v_2\in\QQ$ and $v_1\ne v_2$.)
Theorem~\ref{t:shift} implies, in particular, 
that in $S^2(3) \times  S^2(4)$
the tori $T(a,a+1)$ and $T(a,a+2)$ are Hamiltonian isotopic  for
all sufficiently small $a$, whereas $(a, a+1) \not \simeq (a, a+2)$.
\diam
}
\end{example*}

The paper is organized as follows.
In Section~\ref{s:inv}, we describe the invariants that are used in the
proof of Theorems~\ref{t:t1+} and~\ref{t:tame}, and derive Theorem~\ref{t:tame}.
In Section~\ref{s:t1+} we proof a version of Theorem~\ref{t:tame} for generalized Clifford tori
in~$\CP^n$, and use it to prove Theorem~\ref{t:t1+}.
In Section~\ref{s:con} we 
%F provide the constructions for the proof of Theorem~\ref{t:t1}.
construct Hamiltonian isotopies that provide a proof of Theorem~\ref{t:t1}.
In Sections~\ref{s:Liouville} and~\ref{s:shift}, we prove finer versions of 
Theorems~\ref{t:Liouville} and~\ref{t:shift}, respectively.
Appendix~\ref{a:annulus} provides a refinement of Lelong's inequality for the area of 
holomorphic curves passing through the centre of a ball, that we use in Section~\ref{s:inv}.
Appendix~\ref{a:path} proves an algebraic result used in~Section~\ref{s:con}.

\b
\ni
{\bf Acknowledgments.}
The first draft of this paper was written in Spring~2005, when the first author visited Max Planck Institute Leipzig
and the second author was a PostDoc at Leipzig University.
The paper was finalized during our stay at FIM of ETH Z\"urich in 2008 and~2010
and during the first author's stay at Universit\'e de Neuch\^atel in 2009 and~2011.
We wish to thank these institutions and in particular Dietmar Salamon and Matthias Schwarz 
for their warm hospitality.
% F A somewhat silly sentence added, that in this standard form was prescribed to all CAST members.
The present work is part of the author's activities within CAST, 
a Research Network Program of the European Science Foundation.

%%%%%%%%%%%%%%%%%%%%%%%%%%%%%%%%%%%%%%%%%%%%%%%%%%%%%%%%%%%%%%%%%%%%%%%%%%%%%%%%%%%%%%%%%%%

\section{Symplectic invariants}  \label{s:inv}

\subsection{Displacement energy and $J$-holomorphic discs}
The first Ekeland--Hofer capacity
was a key tool used in~\cite{Ch-tori} for
the classification of product tori in~$\RR^{2n}\?\?$.
This capacity is defined only for subsets of the standard
symplectic space~$\RR^{2n}$.
We shall work with the displacement energy capacity instead,
which is defined for all symplectic manifolds.
In the process of computing the displacement energy for Lagrangian tori,
we bring $J$-holomorphic discs into play, and it is here
that we need the assumption that $(M,\go)$ be tame.

Consider the set $\ch(M)$ of compactly supported smooth
functions $H \colon [0,1] \times M \to \RR$.
Denote by $\Phi_{\?\?H}$ the time 1 map of the Hamiltonian flow generated by~$H$.
Following \cite{Hofer90}, we define a norm on $\ch$ by
\[
\| H \| \,=\, \int_0^1 \left( \max_{x \in M} H(t,x) - \min_{x \in M} H(t,x) \right) dt ,
\]
and the displacement energy of a compact subset $A \subset M$ by
\[
e \left( A,M \right) \,=\, \inf_{H \in \ch}
\Bigl\{ \| H \|
\mid \Phi_H (A) \cap A = \varnothing \Bigr\} ,
\]
assuming $\inf(\varnothing)=\infty$.

Assume that $(M,\go)$ is tame.
Denote by $D$ the closed unit disc in the complex plane~$\CC$,
and by $\cj=\cj(M,\go)$ the set of almost complex structures~$J$ on
$M$ for which there exists a complete Riemannian metric $g$ such
that $J$ and $g$ satisfy (T1) and~(T2).
Let $L$ be a closed Lagrangian submanifold  of~$(M,\go)$.
Given $J \in \cj$, we define $\gs (L,M;J)$ to be the minimal symplectic area $\int_D u^*\go$
of a non-constant $J$-holomorphic map
$u \colon (D, \pp D) \to (M,L)$ if such maps exist, and set $\gs (L,M;J)=\infty$ otherwise.
Since $(M,\go)$ is tame, Gromov's compactness theorem implies
that the minimal area is indeed attained and thus positive~\cite{MSa}.
Define
\[
\gs\left (L,M\right) \,=\, \sup_{J \in \cj} \.\gs (L,M;J) ,
\]
allowing $\gs (L,M)$ to be infinite as well.
It was proved in \cite{Ch-energy} that
\begin{equation}   \label{e:energy}
   \gs \left(L,M\right) \le e \left(L,M\right).
\end{equation}

Recall that $\ua = \min_{1\le i \le n} (a_i)$,
$ \|\ba\| = \ua + \sum_{i=1}^n \?a_i$.

\begin{proposition}  \label{p:e=}
If $(M,\go)$ is tame and $\|\ba\| \.\le b_\gf$, then
$\,e \left( T_\gf (\ba),M  \right) = \ua$.
\end{proposition}

\proof
First we prove that $\,e \left( T_\gf (\ba),M  \right) \le \ua$.
We can assume that $a_1 = \ua$.
We write  $D(\ba)$ for the polydisc $B^2(a_1) \times \dots \times B^2(a_n)$.
%F line added
Let~$U$ be a neighbourhood of $B^{2n}(b)$ such that $\gf \colon U \to M$
is well defined.
Choose $\gve >0$ such that $B^{2n}\bigl( \| \ba +n\gve \|\bigr) \subset U$.
%For each $\gve >0$, 
The torus~$T(\ba)$ can be displaced
from itself by the time 1 flow map of a Hamiltonian function
$H \in \ch\left(D \left( 2 a_1 +\gve,
a_2+\gve, \ldots,  a_n+\gve \right) \right) $
such that $\| H \| < \ua +\gve$,
see e.g.~\cite[p.\hskip.1em 171]{HZ}.
The polydisc $D \left( 2 a_1 +\gve, a_2+\gve, \ldots,  a_n+\gve\right)$
is contained in the ball $B^{2n} ( \.\|\ba\| + n\.\eps)$ and hence in~$U$.
%For $\eps$ small enough, we have $B^{2n} ( \.\|\ba\| + n\.\eps) \subset B^{2n} ( b_\gf)$.
Transferring  $H$ to $(M,\go)$ by means of the chart~$\gf$,
we obtain a Hamiltonian function $H^\gf \!\in\ch(M)$
such that $\| H^\gf \| < \ua +\gve$ and the time 1 flow generated by~$H^\gf$
disjoins the torus $T_\gf(\ba)$ from itself.
Since $\gve$
can be chosen arbitrarily small, it follows that
$\,e \left( T_\gf (\ba),M  \right) \le \ua$.

Denote by $J_0$ the standard complex structure on~$\CC^{n}$.

\begin{lemma}  \label{l:sigma}
Let $L$ be a closed Lagrangian submanifold in~$B^{2n}(b_-)\subset\CC^{n}$,
and let $\gf$ be a symplectic chart such that~$b_\gf \? > b_-$.
Then
\[
\gs \bigl( \gf (L),M \bigr) \.\. \ge
\.\.\min \.\bigl(\gs \.( L,\CC^{n}\?\?; J_0),    \.b_\gf \? - b_-  \bigr).
\]
\end{lemma}

\proof
It suffices to find an almost complex complex structure $J \in \cj$
such that the symplectic area of each non-constant \text{$J$-holomorphic} map
$u \colon (D, \pp D) \to \left(M, \gf(L) \right)$ is at least
$\.\min \.\bigl(\gs \.( L,\CC^{n}\?\? ; J_0),    \.b_\gf \? - b_-  \bigr)$.
We construct such a~$J$ as follows.
Transferring the  complex structure $J_0$ by means of the chart~$\gf$,
we obtain a complex structure $J_0^{\gf}$ on~$B_\gf$.
We claim that $J_0^{\gf}$ extends to an almost complex structure $J \in \cj$ on~$M$.
Indeed, pick an arbitrary~$J_1 \in \cj$.
For each $x\in M$, the space of complex structures $J_{(x)}$ on
the tangent space $T_x M$ satisfying $\go (\xi, J_{(x)} \xi) >0$ for all
$\xi \in  T_x M\setminus\{0\}$ is contractible~\cite{MS}.
Thus there is an almost complex structure $J$ on~$M$ that coincides
with $J_0^{\gf}$ on~$B_\gf$, and with $J_1$ outside a relatively compact
neighbourhood of~$B_\gf$. Then~$J \in \cj$.

Let $u \colon (D, \pp D) \to \left(M, \gf(L) \right)$
be a non-constant \text{$J$-holomorphic} map.
If the image of $u$ is contained in~$B_\gf $, then
$u_\gf\?=\gf^{-1} \circ u \colon (D, \pp D) \to \left( \CC^n, L \right)$
is a non-constant holomorphic map.
Hence $\int_D u^* \omega = \int_D u_\gf^* \.\omega_n \ge \gs \.( L,\CC^{n}\?\? ; J_0)$.

If the image of $u$ is not contained in~$B_\gf $, then
the set $V_\gf=\gf^{-1}(u(D))$ is a real analytic subvariety in
$B(b_\gf)$ intersecting the sphere~$\partial B(b_-)$.
Applying Theorem~\ref{t::a} from Appendix~A
(with $b_-=\pi r^2_-$, $b_\gf=\pi r^2_+$),
we infer that the Riemannian area of $V_\gf$ is at least~$b_\gf \? - b_-$.
Since the Riemannian area of a holomorphic curve in $\CC^{n}$
equals its symplectic area, and the symplectic area of $u$
is not less than the symplectic area of~$V_\gf$,
it follows that the symplectic area of $u$ 
is at least~$b_\gf \? - b_-$.
\proofend

We claim that $\gs \.( T(\ba),\CC^{n}\?\? ; J_0) \ge \ua$.
Let  $u \colon (D, \pp D) \to \left(\CC^{n}\?\?, T(\ba) \right)$ be
a non-constant holomorphic map.
Write $u=(u_1,\ldots,u_n)$, where each
$u_j \colon (D, \pp D) \to \left(\CC, T(a_j) \right)$
is a holomorphic map.
The symplectic area of~$u$ is positive, and it equals the sum of
the symplectic areas of the maps~$u_j$.
Since the symplectic area of~$u_j$ is a non-negative integer
multiple of~$a_j$, the symplectic area of~$u$ is at least~$\ua$.
The torus $T(\ba)$ is contained in the ball $B^{2n} ( |\ba| )$.
By Lemma~\ref{l:sigma},  
$\,\gs\left( T_\gf (\ba),M  \right) \ge \| \ba \| - |\ba| = \ua$.
In view of~\eqref{e:energy},  we conclude that 
$\,e\left( T_\gf (\ba),M  \right) \ge \ua$.
This completes the proof of Proposition~\ref{p:e=}.
\proofend

\subsection{Deformations}
Let $(M,\go)$ be a symplectic manifold.
Denote by $\cl$ the space of closed embedded Lagrangian submanifolds
in $(M,\go)$ endowed with the $C^\infty$-topology.
Given a $\Ham(M,\go)$-invariant function $f$ on $\cl$ taking values in a set~$X$,
we associate with each  $L \in \cl$ a function germ
$ S^f_L  \colon H^1(L;\RR)\to X $ at the point $0\in H^1(L;\RR)$ following~\cite{Ch-tori}.
This construction provides additional invariants of Lagrangian submanifolds.
We use it to prove Theorem~\ref{t:tame}~(ii).
%
%\begin{proposition}  \label{p:ma}
%Let $T_\gf (\ba)$ and $T_\gf (\ba')$
%be product Lagrangian tori
%in a tame symplectic manifold$~(M,\go)$.
%If  $T_\gf (\ba) \sim T_\gf (\ba')$
%and $\,\max \left\{ \|\ba\|, \|\ba'\| \right\} < b_\gf$,
%then $m(\ba) = m(\ba')$.
%\end{proposition}

By Weinstein's Lagrangian Neighbourhood Theorem, there is a symplectomorphism~$g$
from a neighbourhood of $L$ in $M$ to
a fibrewise convex neighbourhood of the zero section of~$T^*\?\?L$
such that the image of $L$ is the zero section~\cite{W-71}.
There is a neighbourhood $V$ of the point $L$ in the space
$\cl$ such that each $L'\in V$ is mapped to the graph
of a closed 1-form $\alpha_{L'}$ on~$L$.
Consider the mapping $w_{L,V}\colon V\to  H^1(L;\RR)$ that sends
$L'\in V$ to the cohomology class of the form $\alpha_{L'}$.
This mapping is locally surjective at~$L$. Denote
by~$w_L$ the germ of $w_{L,V}$ at~$L$.
If two Lagrangian submanifolds $L_0,L_1\in V$ are mapped by $w_{L,V}$ to the same
cohomology class $\zeta\in H^1(L;\RR)$, then they are Hamiltonian isotopic.
Indeed, consider the family of Lagrangian submanifolds $\{L_t\}$
such that $g(L_t)$ is the graph of the 1-form
$\alpha_t=t \alpha_{L_1} + (1-t)\alpha_{L_0}$ for each $t\in \left[0,1\right]$.
Since $[\alpha_t]= \zeta$ for all~$t$, the family $\{L_t\}$
is a Hamiltonian isotopy between $L_0$ and~$L_1$.
Therefore, one can define a mapping germ
$S^f_L \colon H^1(L;\RR)\to X $ at the point $0\in H^1(L;\RR)$
by $S^f_L (\zeta)=f(L')$, where $w_{L}(L')=\zeta$.
In order to prove that the definition of $ S^f_L$ does not depend
on the choice of the symplectomorphism~$g$,
it suffices to give a description of the mapping germ $w_L$
that does not use~$g$.
This description  goes as follows: the evaluation of $w_{L}(L')$ on
a 1-homology class $\lambda\in H_1(L;\ZZ)$ equals
$\int_{[0,1]\times S^1} h^*\omega$, where $h\colon [0,1]\times S^1\to M$ is a
smooth map with image in a tubular neighbourhood of~$L$ such that $h(\{0\}\times S^1)$
is a loop in $L$ representing the class  $\lambda$ and  $h(\{1\}\times S^1)\subset L'$.

It immediately  follows from the definition that $S^f_L$ is $\Ham(M,\go)$-invariant
in the following sense: for each $\psi \in \Ham(M,\go)$, we have
\begin{equation}  \label{e:Sfl}
S^f_{\psi(L)} \,=\,  S_L^f \circ \left( \psi {\mid}_{L}\right)^*,
\end{equation}
and if, moreover, $f$ is $\Symp (M,\go)$-invariant,
then~\eqref{e:Sfl} holds for each  $\psi \in \Symp (M,\go)$.
The displacement energy function $e \left( L \right)=e \left(L,M\right)$
takes values in $\left[ 0, \infty \right[ \cup \{\infty \}$
and is $\Symp (M,\go)$-invariant.

\begin{proposition}  \label{p:Se}
Let $L=T_\gf (\ba)$ be a product  Lagrangian torus in a tame symplectic
manifold. Assume that $\| \ba \| \le b_\gf$.
Then
\[
S_L^e (\zeta) \,=\, e(L) + \min \bigl( l_1(\zeta), \dots, l_{m(\ba)}(\zeta) \bigr),
\]
where $l_1, \dots, l_{m(\ba)} $ are independent linear functions on $H^1(L;\RR)$.
\end{proposition}
\proof
Consider the mapping germ $\theta \colon (\RR^n,0)\to (\cl,L)$, $\bs\mapsto T_\gf (\ba+ \bs)$.
The composition $A=w_L\circ\theta \colon (\RR^n,0)\to (H^1(L;\RR),0)$
is a linear isomorphism germ.
%F 1 line added
Choose $\gve >0$ so small that $\gf \colon B^{2n}(b_\gf +\gve) \to M$ is defined.
For $\bs$ small enough,  we have $\|\ba+\bs\| \le b_\gf + \gve$ and hence,
by Proposition~\ref{p:e=},
\begin{equation}
\label{e:e1}
e\left( T_\gf (\ba + \bs) \right) \,=\, \min \bigl ( a_1+s_1, \dots, a_n+s_n \bigr).
\end{equation}
We can assume, after reordering the coordinates, that
\[
\ua = a_1 =  \dots = a_{m(\ba)} < a_{m(\ba)+1} \le \dots \le a_n.
\]
For $\bs$ sufficiently small (say, such that the absolute values of all
its components do not exceed  $\tfrac{1}{2}(a_{m(\ba)+1} - a_{m(\ba)})$),
in view of~\eqref{e:e1} we have
\begin{equation}  \label{e:e2}
e\left( T_\gf (\ba + \bs) \right) \,=\, 
\ua + \min   \bigl ( s_1, \dots, s_{m(\ba)} \bigr)
\,=\, e(L) + \min   \bigl ( \pi_1(\bs), \dots, \pi_{m(\ba)}(\bs) \bigr),
\end{equation}
where $\pi_i \colon \RR^n \to \RR$ is the projection
onto the $i$-th coordinate axis, $\pi_i(\bs)=s_i$.
Since $S_L^e (\zeta)= e\bigl( T_\gf (\ba + A^{-1}(\zeta)) \bigr)$,
it follows~from~\eqref{e:e2} that
\[
S_L^e (\zeta) \,=\, e(L) + \min \bigl( l_1(\zeta), \dots, l_{m(\ba)}(\zeta) \bigr)
\]
where $l_1=\pi_1\circ A^{-1}, \dots,  l_{m(\ba)}=\pi_{m(\ba)}\circ A^{-1}$
are independent linear functions on~$H^1(L;\RR)$.
\proofend

\proofof{Theorem~\ref{t:tame}~(ii)}
Denote $L=T_\gf(\ba)$, $L'=T_{\gf'}(\ba')$.
It follows from Proposition~\ref{p:e=} and symplectic invariance of
displacement energy that
\[ \ua\,=\,e \left(L,M  \right) \, = \,
e \left( L',M  \right) \,=\, \ua'.\]
According to Proposition~\ref{p:Se},
the cohomology classes $\zeta\in H^1(L;\RR)$ such that $S_L^e(\zeta)=\ua$
form a vector space germ $W$ of dimension~$n-m(\ba)$,
and the cohomology classes $\zeta'\in H^1(L';\RR)$ such that $S_{L'}^e(\zeta')=\ua$
form a vector space germ $W'$ of dimension~$n-m(\ba')$.
If $L'= \psi (L)$ for some $\psi \in \Symp (M,\go)$, then
$S^e_{L'} =  S_L^e \circ A_\psi$, where
$A_\psi= \left( \psi {\mid}_{L}\right)^*$
is a linear isomorphism between $H^1(L;\RR)$ and $H^1(L';\RR)$.
Hence $A_\psi (W) = W'$, and therefore $m(\ba)=m(\ba')$.
\proofend

\subsection{Symplectic area class and Maslov class}
Given a Lagrangian submanifold~$L$ of a symplectic manifold~$(M,\go)$,
one can consider two relative cohomology classes: 
the symplectic area class
$\gs_L \in H^2(M,L;\RR)$ represented by the 2-form~$\omega$,
and the Maslov class $\mu_L\in H^2(M,L;\ZZ)$, defined as in~\cite{V-87}.
Both $\gs$ and $\mu$ are symplectically invariant
in the sense that $\gs_{\psi (L)} = \psi^*\gs_L$ and $\mu_{\psi (L)} = \psi^* \mu_L$
for each symplectomorphism~$\psi$.
These classes define homomorphisms from $\pi_2(M,L)$ to $\RR$ that
we shall also denote by $\gs_L$ and~$\mu_L$.
Define the subgroup $\Gamma (L)\subset \RR$ to be
the image of the subgroup $\ker (\mu_L) \subset \pi_2(M,L)$
under the homomorphism $\gs_L \colon \pi_2(M,L) \to \RR$.
Since $\gs_L$ and $\mu_L$ are symplectically invariant, so is~$\Gamma (L)$:
\begin{lemma}  \label{l:gammaL}
Let $L,L'$ be Lagrangian submanifolds of ~$(M,\go)$.
If $L\sim L'$, then $\Gamma (L)=\Gamma (L')$.
\end{lemma}

Theorem~\ref{t:tame}~(i) is a corollary of
Lemma~\ref{l:gammaL} and the following assertion:

\begin{lemma}  \label{l:gamma}
Let $T_{\gf}(\ba)$ be a product Lagrangian torus in a symplectically
aspherical symplectic manifold~$(M,\go)$.
Then $\Gamma \bigl(T_{\gf}(\ba)\bigr)=\Gamma (\ba)$.
\end{lemma}

\proof
For $i \in \{ 1, \dots, n \}$, let $D_i$ be a disc in $\RR^{2n}$
with boundary on $T(\ba)$ such that the projection of~$D^i$
to the $i$-th factor in $\RR^{2}\times \cdots \times\RR^{2}=\RR^{2n}$
is the  disc in $\RR^2$ bounded by the circle $T(a_i)$,
and the projections to other factors are points.
Denote by $\hat D_i$ the element of
$\pi_2 \left( \RR^{2n},T(\ba) \right)$ represented by~$D_i$.
The classes $\hat D_1,\ldots,\hat D_n$ generate
the free Abelian group $\pi_2 \left( \RR^{2n},T(\ba) \right)$.
Denote $\tilde D_i=\gf_*\hat D_i \in \pi_2 ( M,L )$ where $L := T_\gf(\ba)$.
For each~$i$, we have $\gs_{\.T(\ba)}\?(\hat D_i)=a_i$, $\mu_{\.T(\ba)}\?(\hat D_i) = 2$,
and hence $\gs_{L}(\tilde D_i)=a_i$,  $\mu_{L}(\tilde D_i) = 2$.
The group $\pi_2 \left( M,L \right)$ is the direct sum of $\pi_2(M)$
and the subgroup generated by the elements~$\tilde D_i$.
Since $(M,\go)$ is symplectically aspherical and $\mu_L|^{\vphantom1}_{H_2(M;\ZZ)}  = 2\.c_1 (\go)$ 
(see~\cite{V-87}), 
both $\gs_L$ and $\mu_L$ vanish on~$\pi_2(M)$.
The kernel of $\mu_L$ is the direct sum of $\pi_2(M)$
and the subgroup generated by the differences~$\tilde D_i-\tilde D_j$,
where $i,j \in \{1,\ldots,n\}$ and $j$ is such that $\ua=a_j$.
Therefore, $\gs_L (\ker \mu_L)$ consists of all integer combinations
of the numbers $a_i-\ua=\gs(\tilde D_i-\tilde D_j)$.
\proofend
%

%%%%%%%%%%%%%%%%%%%%%%%%%%%%%%%%%%%%%%%%%%%%%%%%%%%%%%%%%%%%%%%%%%%%%%%%%%%%%%%%%%%

\section{Proof of Theorem~\ref{t:t1+}}  \label{s:t1+}

\subsection{Generalized Clifford tori in $\CP^n$}  \label{s:clifford}

We consider a certain class of product Lagrangian tori in
the complex projective space, the so-called {\it generalized Clifford tori}.
Identify the symplectic space $(\RR^{2n},\go_n)$ with~$\CC^{n}\?$,
the complex coordinates being $z_1=x_1\?+i\.y_1, \ldots, z_n = x_n\?+i\.y_n$.
Consider the diagonal action of the Lie group $U(1)$ on the space~$\CC^{n}\?$.
For each $b>0$, the sphere $S^{2n-\?1\?}(b)=\partial B^{2n}(b)$ is
invariant under this action.
Denote by $\CC P^{\.n-\?1\?}(b)$ the quotient $S^{2n-\?1\?}(b)/U(1)$.
%The orbits contained in $S^{2n-\?1\?}(b)$ are tangent to the kernel of
%the $2$-form~$\go_n|_{S^{2n-\?1\?}(b)}$.
The restriction of the symplectic form $\go_n$ to ${S^{2n-\?1\?}(b)}$
is the pullback of a certain symplectic form
$ \go_{n-\?1}^b$ on~$\CC P^{\.n-\?1\?}(b)$.
This form is a multiple of the Fubini--Study form.

If $\ba\in\RR_+^{n}$ and $|\ba|=b$, then the torus $T(\ba)$
is contained in the sphere~$ S^{2n-\?1\?}(b)$.
Moreover, $T(\ba)$ is invariant under the action of~$U(1)$.
Therefore, the quotient $\T(\ba)=T(\ba)/U(1)$ is a
Lagrangian $(n-\?1)$-torus in~$\CC P^{\.n-\?1\?}(b)$.
It is called a generalized Clifford torus.

Denote by $Z_n(b)$ the complex hypersurface
$$
\bigr(S^{2n-\?1\?}(b)\cap\{z_n\?=\?0\}\bigl)/U(1) \cong \CC P^{\.n-\?2\?}
$$
in~$\,\CC P^{\.n-\?1\?}(b)$,
and by  $\Bcirc^{2n-2}(b)$ the open ball~$\mathrm{Int} (B^{2n-2}(b))$.
The tori $\T(\ba)$ are product tori:

\begin{proposition}  \label{p:cpn}
There is a symplectomorphism
$$ 
\gf^b_{n-\?1} \colon 
\bigl(\Bcirc^{2n-2}(b),\go_{n-\?1}\bigr)
\,\to\,
\bigl(\CC P^{\.n-\?1\?}(b)\setminus Z_n(b),
  \go_{n-\?1}^b\bigr)
$$
that maps each product torus $T(a_1,\dots,a_{n-\?1})$
contained in $\Bcirc^{2n-2}(b)$ to the torus
$\T(a_1,\dots,a_n)$, 
where $a_n= b - a_1- \dots - a_{n-\?1}$.
\end{proposition}

\proof
Denote by $W$ the subset of $S^{2n-\?1\?}(b)$ formed
by points with $z_n$ coordinate positive real.
Consider the projection of $ \CC^n$ onto $\CC^{\.n-\?1\?}$
defined by forgetting the last coordinate.
Restricting this projection to $W$ be obtain a 
diffeomorphism $\psi \colon W \to \Bcirc^{2n-2}\?\?(b)$.
We claim that $\psi$ is a symplectomorphism from 
$(W, \go_n|^{\vphantom1}_W)$ onto 
$\bigl(\Bcirc^{2n-2}\?\?(b),\go_{n-\?1}\bigr)$.
This statement is equivalent to the assertion that
the restriction of the 2-form $dx_n\wedge dy_n$ to $W$ vanishes.
The latter follows since $y_n$ vanishes on~$W$.

The manifold $S^{2n-\?1\?}(b)\setminus \{z_n\?=\?0\}$
is foliated by the orbits of the $U(1)$-action.
Each of these orbits intersects $W$ exactly once,
and the intersection is transverse.
Therefore, symplectic reduction gives rise to a canonical 
symplectomorphism $\psi'$ from 
$(W, \go_n|^{\vphantom1}_W)$ onto
$\bigl(\CC P^{\.n-\?1\?}(b)\setminus Z_n(b), \go_{n-\?1}^b\bigr)$.

The composition $\gf^b_{n-\?1\?}=\psi'\circ\psi^{-1}$ 
is the required symplectomorphism. 
%In coordinates, it is given by
%???
%\[
%\gf_{n-\?1\?}^b \colon (z_1, \dots, z_{n-\?1}) \mapsto \left[ z_1 : \dots : z_{n-\?1} :
%\sqrt{\frac b\pi -|z_1|^2 - \dots - |z_{n-\?1}|^2} \right] .
%\]
%
To prove the assertion concerning Lagrangian tori,
it suffices to observe that the image of
$T(a_1,\dots,a_{n-\?1})$ under the symplectomorphism $\psi^{-1}$ is 
the torus $T(a_1,\dots,a_{n-\?1})\times \sqrt{a_n/\pi}$,
and that the $U(1)$-orbits passing through the latter torus 
form the torus $T(a_1,\dots,a_n)$.
\proofend
%

%Y Remark eliminated

\begin{proposition}  \label{p:cpn2}
Let $\ba,\.\ba'\in\RR_+^{n}$ be such that $|\ba|=|\ba'|$.
Consider the Lagrangian tori $\,\T(\ba)$, $\T(\ba')$ in $\,\CC P^{\.n-\?1\?}(|\ba|)$.
If $\,\T(\ba)\sim \T(\ba')$, then $\ba\simeq\ba'$.
\end{proposition}

\proof
By Theorem~\ref{t:tame}~(ii) we have $\ua = \ua'$ and $m(\ba)=m(\ba')$.
In view of Lemma~\ref{l:gammaL},
it remains to show that $\Gamma(\T(\ba))=\Gamma(\ba)$.
Let $\hat D_1,\ldots,\hat D_{n-\?1}$
be the elements of the group
$\pi_2 \left( \RR^{2n-2},T(a_1,\ldots,a_{n-\?1}) \right)$
defined as in the proof of Lemma~\ref{l:gamma}.
%F The inverse of the symplectomorphism $\psi^b_n$   
The symplectomorphism $\gf^{|\ba|}_{n-1}$ 
sends these classes
to the classes $\tilde D_1,\ldots,\tilde D_{n-\?1}$
in $\pi_2(\CC P^{\.n-\?1\?}(|\ba|), \T(\ba))$.
For each~$i$, we have $\gs_{\t(\ba)}\?(\tilde D_i)=a_i$,
$\mu_{\t(\ba)}\?(\tilde D_i) = 2$.
The free Abelian group $\pi_2(\CC P^{\.n-\?1\?}(|\ba|), \T(\ba))$
is generated by the classes $\tilde D_1,\ldots,\tilde D_{n-\?1}$,
and the class $[\CC P^1]$ represented by a complex line in
the complex projective space.

We have $\mu_{\t(\ba)}\?\?\bigl(\.[\CC P^1]\.\bigr) = 2n$,
since the value of the Maslov class on $\CC P^1$
is twice the value of~$c_1(T\.\CC P^{\.n-\?1\?})$.
We claim that $\gs_{\t(\ba)}\?\?\bigl(\.[\CC P^1]\.\bigr)=|\ba|$.
Indeed, let $\CC P^1 \subset \CP^{n-1}$ be the quotient
of the sphere $\{z_2=\dots=z_{n-\?1}=0\} \cap 
%F \CC P^{\.n-\?1\?}(|\ba|)$
S^{2n-1}(|\ba|)$
by the diagonal action of~$U(1)$.
The symplectomorphism 
%F $\psi^{\?|\ba|}_n$
$\gf_{n-1}^{|\ba|}$
identifies the complement of a 
%F circle 
point
in $\CC P^1$ with the open symplectic disc
$\Bcirc^{2n-2}\?\?(|\ba|)\cap \{z_2=\dots=z_{n-\?1}=0\}$.
This disc has area~$|\ba|$, and
hence the integral of the symplectic form
$ \go_{n-\?1}^{\?|\ba|}$ over $\CC P^1$ equals~$|\ba|$.

Define $\tilde D_n=[\CC P^1]-\sum_{i=1}^{n-\?1} \tilde D_i$.
The group  $\pi_2(\CC P^{\.n-\?1\?}(|\ba|), \T(\ba))$
is generated by the classes $\tilde D_1,\ldots,\tilde D_n$,
and we have
$\gs_{\t(\ba)}\?(\tilde D_n)=a_n$,
$\mu_{\t(\ba)}\?(\tilde D_n) = 2$.
The kernel of $\mu_{\t(\ba)}$ is  generated by the
differences~$\tilde D_i-\tilde D_j$,
where $i,j\in\{1,\ldots,n\}$ and $j$ is such that $\ua=a_j$.
Therefore, $\gs_{\t(\ba)} (\ker \mu_{\t(\ba)})$
consists of all integer combinations
of the numbers $a_i-\ua=\gs(\tilde D_i-\tilde D_j)$.
\proofend

%F proof added
\subsection{Proof of Theorem~\ref{t:t1+}}
Arguing by contradiction, we suppose that $T(\ba) \approx T(\ba')$ in $B^{2n}(b)$.
By Theorem~\ref{t:tame}~(ii), with $(M,\go)$ a large ball and $\gf, \gf'$ the identity embeddings,
we have $\ua = \ua'$.
We can assume that $\| \ba \| \ge \| \ba' \|$.
Since $\ua = \ua'$ and, by hypothesis, $| \ba | \neq | \ba' |$,
we have $\| \ba \| - \| \ba' \| = |\ba | - |\ba' | >0$.
By hypothesis we have $|\ba| \le b < \| \ba \|$.
Thus $|\ba' | < |\ba| \le b < |\ba|+\ua$.
Choose $c'<c$ such that
$$
b < c' < c < |\ba|+\ua .
$$
Define $a_{n+1} := c-|\ba|$ and $a_{n+1}' := c-|\ba'|$.
Then $a_{n+1} < a_{n+1}'$ and $a_{n+1} = c-|\ba| < \ua$.
Therefore,
\begin{equation} \label{e:ineq}
\min \{ a_1, \dots, a_n, a_{n+1} \} \,=\, 
%\min \{ \ua, a_{n+1}' \} \,=\,
a_{n+1} < \min \{ \ua, a_{n+1}' \} \,=\,
\min \{ a_1', \dots, a_n', a_{n+1}' \} .
\end{equation}
Recall that $T(\ba) \approx T(\ba')$ in $B^{2n}(b)$.
Cutting off the Hamiltonian function that generates this isotopy,
we construct a Hamiltonian isotopy supported in $B^{2n}(c')$ that moves 
$T(\ba)$ to~$T(\ba')$.
The symplectomorphism $\gf_{n}^c$ from Proposition~\ref{p:cpn}
transfers this isotopy to a Hamiltonian isotopy of $\CP^n(c)$.
It moves $\T(a_1, \dots, a_n,a_{n+1})$ to $\T(a_1', \dots, a_n',a_{n+1}')$.
By Proposition~\ref{p:cpn2}, 
$\min \{ a_1, \dots, a_n, a_{n+1} \} = \min \{ a_1', \dots, a_n', a_{n+1}' \}$,
in contradiction to~\eqref{e:ineq}.
\proofend

\section{Constructions of Hamiltonian isotopies}  \label{s:con}

\subsection{Proof of Theorem~\ref{t:t1}~(i)}
The unitary group $U(n)$ acts on $\CC^n$ preserving the symplectic form~$\omega_n$.
Since a permutation of coordinates $z_1,\ldots, z_n$ is a unitary map and
the group $U(n)$ is path-connected, there is a smooth family $\{\Phi_t\}$,
$t\in[0,1],$ of unitary maps such that $\Phi_0={\tt id}$ and
$\Phi_1(T(\ba))=T(\ba')$. The flow $\{\Phi_t\}$ is Hamiltonian because $\CC^n$
is simply-connected.
\proofend

\subsection{}
The proof of Theorem~\ref{t:t1}~(ii) relies on
the following lemma, which represents a special case of Theorem~\ref{t:t1}~(ii).
\begin{lemma}  \label{l:step1}
For each positive $a$, $c$, and $d$,
the tori $T(a,a+c,a+d)$ and $T(a,a+c+d,a+d)$ are
Hamiltonian isotopic in the ball $B^6 (4a+ c+2d)$.
\end{lemma}

\proof
Let $\. W= \bigl\{ (z_1,z_2) \?\?\in \?\CC^2 \?\bigm| |z_1|\?<\?|z_2|\.\.\bigr\} $.
Consider the map
\[
\Psi \colon W \to \CC^2, \quad (z_1,z_2) \mapsto
\left(\frac{z_1 z_2}{|z_2|} ,\frac{z_2 \sqrt{|z_2|^2 \?\?-|z_1|^2} }{|z_2|}\right).
\]
It is injective, and
its image is the complement of the complex line $\{\.z_2=0\.\}$.
We claim that $\Psi$ preserves the symplectic form
$\omega_2= dx_1\wedge dy_1 + dx_2\wedge dy_2$.
Indeed, write $z_1= e^{2\pi i\.\.\theta_1\!}\sqrt{\rho_1/\pi} $, 
$z_2= e^{2\pi i\.\.\theta_2\!}\sqrt{\rho_2/\pi} $,
with $\theta_1,\theta_2$ in $S^1=\RR/\ZZ$ and $\rho_1,\rho_2$ nonnegative real.
For nonzero values of $z_2$, we have
$\omega_2= d\rho_1\wedge d \theta_1 + d\rho_2\wedge d \theta_2$ and
\[
\Psi (\rho_1, \theta_1, \rho_2,  \theta_2) =
(\rho_1, \theta_1+\theta_2, \rho_2- \rho_1,  \theta_2).
\]
Clearly, $\Psi$ is symplectic outside the complex line $\{\.z_2=0\.\}$, and hence,
by continuity, on the whole of~$W$.
A product torus $T(a_0, a_0+b_0)\subset W$ is mapped by
$\Psi$ to the torus $T(a_0, b_0)$.

The torus $T(a,a+c+d,a+d)$ is Hamiltonian isotopic, through a unitary isotopy,
to the torus $T(a+d,a+c+d,a)$ in the ball $B^6(3a+ c+2d)$.
Therefore, it suffices to prove that the tori $T(a,a+c,a+d)$ and $T(a+d,a+c+d,a)$
are Hamiltonian isotopic in~$B^6 (4a+ c+2d)$.

Consider the map
$\Psi_{\!+}=\Psi\times {\tt id} \colon W\?\times\?\CC\to\CC^3$.
We have $\Psi_{\!+}(T(a,a+c,a+d))= T(a,c,a+d)$ and
$\Psi_{\!+}(T(a+d,a+c+d,a))= T(a+d,c,a)$.
The Hamiltonian function $H=\frac{\pi}{2}(x_1 y_3 -x_3 y_1)$
gives rise to a unitary Hamiltonian flow $\{\Phi_t\}$
that does not change the complex coordinate $z_2$.
We have $\Phi_1(z_1,z_2,z_3)=(z_3,z_2,-z_1)$.
In particular,  $\Phi_1$ maps $T(a,c,a+d)$ to $T(a+d,c,a)$.
Multiplying $H$ by an appropriate cutoff function, we construct a Hamiltonian~$H'$,
 compactly supported in $\CC^3\setminus \{\.z_2=0\.\}$, whose
flow $\{\Phi'_t\}$  moves the torus $ T(a,c,a+d)$ in exactly the same way as the flow~$\{\Phi_t\}$.
Consider the Hamiltonian flow $\{\Phi^+_t\}$ on $\CC^3$ generated by the
Hamiltonian function $H'\circ\Psi_{\!+}$.
%
%F line added
This flow is compactly supported in $W \times \CC$, where
$\Phi_t^+ = \Psi_+^{-1} \circ \Phi_t' \circ \Psi_+$.
In particular, 
$\Phi^+_t(T(a,a+c,a+d))= \Psi_{\!+}^{{-}1}(\Phi_t(T(a,c,a+d)))$ for all values of~$t$,
and $\Phi^+_1(T(a,a+c,a+d))= T(a+d,a+c+d,a)$.
It remains to show that each torus $\Phi^+_t(T(a,a+c,a+d))$ is contained in $B^6 (4a+ c+2d)$.

Let $(z_1,z_2,z_3)\in \Phi^+_t(T(a,a+c,a+d))$.
We are to prove that $\pi\.(\.|z_1|^2\?+|z_2|^2\?+|z_3|^2)\le 4a+c+2d$.
The point $\Psi_{\!+}(z_1,z_2,z_3)= (z'_1,z'_2,z_3)$ belongs to the
torus~$\Phi_t(T(a,c,a+d))$.
Since $T(a,c,a+d)$ is contained in the sphere $\partial B^6 (2a+ c+d)$ and $\Phi_t$ is unitary,
it follows that $(z'_1,z'_2,z_3)\in \partial B^6 (2a+c+d)$.
Hence $\pi\.(\.|z'_1|^2\?+|z'_2|^2\?+|z_3|^2)= 2a+c+d$.
By the construction of $\Phi_t$, we have $\pi|z'_2|^2=c$.
The definition of the map $\Psi_{\!+}$
implies that $|z'_1| =|z_1|$ and $|z_2|^2= |z'_1|^2\?+ |z'_2|^2$.
Therefore,
\begin{eqnarray*}
\pi\.(\.|z_1|^2\?+|z_2|^2\?+|z_3|^2) &=& 2a+c+d+ \pi |z'_1|^2 \\
&=&
4a+2c+2d-\pi|z'_2|^2 - \pi|z_3|^2 \\
&=& 
4a+c+2d - \pi|z_3|^2 \,\le\, 4a+c+2d ,
\end{eqnarray*}
%F 5 words added
as we wished to show.
\proofend

\begin{lemma}  \label{l:step2}
Let $\bc=(c_1,\dots,c_k)$ and $\bc'=(c'_1,\dots,c'_k)$ be vectors in~$\RR_+^k$, $k \ge 2$, 
such that, for some different indices $i,j\in \{1,\ldots,k\}$,
we have $c_i'=c_i+c_j$,   and $c_l'=c_l$ for $l\ne i$.
For each $n>k$ and each positive~$a$,
the $n$-dimensional tori
$T(\bp)=T(a,\dots,a, a+c_1,\ldots, a+c_k)$ and
$T(\bp')=T(a,\dots,a, a+c'_1,\ldots, a+c'_k)$
are Hamiltonian isotopic in the ball $B^{2n} (\|\bp'\|)$.
\end{lemma}
\proof
We may assume that $i=1$ and $j=2$
after applying to the tori $T(\bp)$ and $T(\bp')$ unitary isotopies
that swap the  complex coordinates $z_{n-k+1}$ and~$z_{n-k+i}$,
$z_{n-k+2}$ and~$z_{n-k+j}$.
By Lemma~\ref{l:step1}, there is a Hamiltonian isotopy on $\CC^3$ that moves
the torus $L_0=T(a,a+c_1,a+c_2)$ to $L_1=T(a,a+c_1+c_2,a+c_2)$
through tori $L_t$ belonging to $B^6(4a+c_1+2 c_2)$.
The tori $L'_0=T(\bp)$ and $L'_1=T(\bp')$ are Hamiltonian isotopic
through the family $L'_t= T(a,\ldots,a)\times L_t \times T(a+c_3,\ldots,a+c_k)$.
All the tori $L'_t$ are contained in the ball
$$
B^{2n}\bigl((n+1)\.a+|\bc|+c_2 \bigr)  = B^{2n} \bigl(\|\bp'\|\bigr).
$$
\proofend

\subsection{Proof of Theorem~\ref{t:t1}~(ii)}
After applying appropriate unitary isotopies to the tori $T(\ba)$ and $T(\ba')$,
we may assume that the first
$m(\ba)$ components of both $\ba$ and $\ba'$ equal~$\ua$.
Let $k=n-m(\ba)$.
Write
\[
 T(\ba)=T(\ua,\dots,\ua, \ua+d_1,\ldots, \ua+d_k),\quad
 T(\ba')=T(\ua,\dots,\ua, \ua+e_1,\ldots, \ua+e_k),
\]
where  $\bd=(d_1,\dots,d_k)$ and
$\be=(e_1,\ldots,e_k)$ are vectors in~$\RR_+^k$.
If $k$ equals~$1$, then the hypothesis
$\Gamma(\ba)=\Gamma(\ba')$ implies that $\ba = \ba'$,
and there is nothing to prove.
Assume that $k\ge2$.

Suppose that there is a sequence
$\bd = \bd^0, \bd^1, \bd^2, \ldots, \bd^{\ell} = \be$
of vectors in $\RR_+^k$ with the following property:
 for each $s\in \{1,\dots,\ell\}$,
the vector $\bd^s=(d^s_1,\ldots,d^s_k)$ is obtained from the vector $\bd^{s-1}$
either by swapping two of the components, or by
adding to the $i$-th component the $j$-th component, or by
subtracting from the $i$-th component the $j$-th component.
For $s\in\{0,\ldots,\ell\}$, define $\ba^s=(\ua,\dots,\ua, \ua+d^s_1,\ldots, \ua+d^s_k)$.
Consider the  sequence of tori
$ T(\ba)=T(\ba^0),T(\ba^1),\ldots,T(\ba^\ell)= T(\ba')$.
For each $s\in\{1,\ldots,\ell\}$, the tori $T(\ba^{s-1})$ and $T(\ba^{s})$
are Hamiltonian isotopic inside the ball $B^{2n} \bigl(\max(\|\ba^{s-1}\|,\|\ba^s\|)\bigr)$.
Indeed, if $\bd^{s-1}$ and $\bd^{s}$ are related by a swap of components, then there is
a unitary isotopy; otherwise, we
apply Lemma~\ref{l:step2} with either $\bc=\bd^{s-1}$, $\bc'=\bd^{s}$,
or $\bc=\bd^{s}$, $\bc'=\bd^{s-1}.$

It thus suffices to show that such a sequence $\bd^s$ indeed exists
and that, moreover,
we have $\|\ba^{s}\|\le \|\ba\|$ or $\|\ba^{s}\|\le \|\ba'\|$ for all~$s$.
To this end, we apply Theorem~\ref{t:low} from Appendix~B
(we refer the reader to the definitions therein).
The theorem is applicable because the condition $\Gamma(\ba)=\Gamma(\ba')$
means exactly $\langle \bd \.\rangle = \langle \be \rangle$.
Consider the sequence (path) $\bd = \bd^0, \bd^1, \bd^2, \ldots, \bd^{\ell} = \be$
whose existence is guaranteed by Theorem~\ref{t:low}.
Since it is low, it follows that, for each~$s$,
we have $|\bd^s|\le |\bd|$ or $|\bd^s|\le |\be|$, and
hence $\|\ba^s\|\le \|\ba\|$  or $\|\ba^s\|\le \|\ba'\|$.
Thus, this sequence has the required properties and
the proof of Theorem~\ref{t:t1} is complete.
\proofend

\section{Spaces of symplectic charts and product tori}
\label{s:Liouville}

Given $b>0$, denote by $\Emb \left( B^{2n}(b),M,\go \right)$
the space of symplectic charts
$\gf \colon B^{2n}(b) \to (M,\go)$, endowed with the $C^\infty$-topology.
By Darboux's theorem, this space is nonempty
at least for sufficiently small~$b$.
The {\it Gromov radius}\/ $\rho (M,\go)$ of $(M,\go)$ is defined as
the supremum of all $b$ such that $\Emb \left( B^{2n}(b),M,\go \right)$
is nonempty (we allow $\rho (M,\go)=\infty$).
For computations and estimates of $\rho \.(M,\go)$
we refer to~\cite{Sch-book} and the references therein.
It has been conjectured that the space $\Emb \left( B^{2n}(b),M,\go \right)$
is connected for all closed symplectic manifolds and all~$b>0$.
This has been proved for certain closed $4$-manifolds
and also for the symplectic $4$-ball~$\Bcirc(c)$,
see~\cite{McDuff-96}.

\begin{theorem}  \label{t:rho}
Let $T_\gf(\ba)$ and $T_{\gf'}(\ba')$
be two Lagrangian product tori
in a  symplectically aspherical tame symplectic manifold~$(M,\go)$.

\s
{\rm (i)}
Let $b_- = \min \left\{ b_\gf, b_{\gf'} \right\}$ and
$b_+ = \max \left\{ b_\gf, b_{\gf'} \right\}$.
Assume that the space
$\Emb \left( B^{2n}(b_-),M,\go \right)$ is path-connected and that
$ \max \bigl\{  \| \ba \|, \|\ba'\| \bigr\} \le b_+$.
Then the conditions $\ba \simeq \ba'$, $T_\gf(\ba) \sim T_{\gf'}(\ba')$,
$T_\gf(\ba)\approx T_{\gf'}(\ba')$ are equivalent one to another.

\s
{\rm (ii)}
Assume that the space
$\Emb \left( B^{2n}(b),M,\go \right)$ is path connected
for all values of~$b$
and that $ \max \bigl\{  \| \ba \|, \|\ba'\| \bigr\} < \rho (M,\go)$.
Then the conditions $\ba\simeq\ba'$, $T_\gf(\ba)\sim T_{\gf'}(\ba')$,
$T_\gf(\ba)\approx T_{\gf'}(\ba')$ are equivalent one to another.
\end{theorem}

\proof
First we prove statement~(i).
If $T_\gf(\ba)\approx T_{\gf'}(\ba')$, then
$T_\gf(\ba)\sim T_{\gf'}(\ba')$ by definition.
If $T_\gf(\ba)\sim T_{\gf'}(\ba')$, then $\ba\simeq\ba'$
by Theorem~\ref{t:tame}.
Let~$\ba\simeq\ba'$.
We can assume that $b_- = b_\gf$ and $b_+ = b_{\gf'}$.
It follows from Theorem~\ref{t:t1} that
$T_{\gf'}(\ba)\approx T_{\gf'}(\ba')$.
Since $\Emb \left( B^{2n}(b_-),M,\go \right)$ is path-connected,
there exists a smooth family $\{\gf_s\}$, $s \in [0,1]$,
of symplectic embeddings $B^{2n}(b_-) \to (M,\go)$ such that
$\gf_0 = \gf$ and $\gf_1$ coincides with $\gf'$ on~${B^{2n}(b_-)}$.
Then there is a Hamiltonian isotopy
$\{\Psi_s\}$, $s \in [0,1]$, of~$(M,\go)$ such that
$\Psi_s\circ\gf=\gf_s$ for all~$s$.
We have $\Psi_1(T_\gf(\ba))=T_{\gf'}(\ba)$.
Therefore, $T_{\gf}(\ba)\approx T_{\gf'}(\ba)\approx T_{\gf'}(\ba')$.

The statement~(ii) will follow from the statement~(i)
if we show that, for each $b,b'$ satisfying
$0<b<b'<\rho (M,\go)$, every symplectic embedding
$\gf \colon B^{2n}(b) \to (M,\go)$
extends to a symplectic embedding
$\colon B^{2n}(b') \to (M,\go)$.
Pick $\gf_+\in\Emb \left( B^{2n}(b'),M,\go \right)$.
Denote by $\gf'_+$ the restriction of  $\gf_+$ to~$\gf \colon B^{2n}(b)$.
Since $\Emb \left( B^{2n}(b),M,\go \right)$ is path-connected,
we conclude, arguing as above, that there is $\Psi\in \Symp(M,\go)$
such that $\Psi\circ \gf = \gf'_+$.
Then $\Psi^{-1}\circ \gf_+\in \Emb \left( B^{2n}(b'),M,\go \right)$ is
an extension of~$\gf$.
\proofend

%\subsection{Proof of Theorem~\ref{t:Liouville}}  \label{ss:liouville}

%
\begin{proposition}  \label{p:emb}
For a forward complete Liouville manifold~$(M,\go)$,
the space $\Emb \left( B^{2n}(b),M,\go \right)$
is nonempty and path-connected for each~$b>0$.
\end{proposition}

\proof
Let $X$ be a forward complete Liouville field on $(M,\go)$.
Denote by $\{f_t\}$, $t \ge 0$, its forward flow.
Assume that the space $\Emb\left( B^{2n}(b),M,\go \right)$
is nonempty and pick $\gf\in\left( B^{2n}(b),M,\go \right)$.
Since $\left( f_t \right)^* \go = e^t \go$ for all $t \ge 0$,
the map
\[  B^{2n}(e^{2t} b)\to M, \quad x\mapsto f_{2t} ( \gf  (e^{-t} x) ) \]
is a symplectic embedding, and hence the space
$\Emb\left( B^{2n}(b_+),M,\go \right)$ is nonempty for all $b_+>b$.

Let $\gf, \gf' \colon B^{2n}(b) \to (M,\go)$.
We prove that $\gf$ and $\gf'$ are homotopic through symplectic embeddings.
After composing $\gf'$ with an appropriate Hamiltonian symplectomorphism
of $(M,\go)$, we can assume that $\gf (0) = \gf' (0)$.
Since each element of the linear symplectic group $\Sp (2n;\RR)$
can be realized as linearization of a Hamiltonian symplectomorphism
preserving the point~$\gf (0)$,
we can also assume that $d \gf (0) = d \psi (0)$.
There is a symplectic isotopy
$\{F_t\}$, $t \in [0,1]$, of $B^{2n}(b)$
such that $F_0 = \id$ and $\psi \circ F_1$ coincides with $\gf$ on
$B^{2n}(b')$ for some $b'\in\left]\.0,b\.\right[$,
see~e.g.\ Appendix~A.1 of \cite{HZ} or the proof of Lemma~2.2 in~\cite{RS}.
Therefore, we may assume that $\gf=\psi$ on~$B^{2n}(b')$.

Consider smooth families $\{\Phi_t\}, \{\Psi_t\}$, $t \ge 0$,
of embeddings $ B^{2n}(b) \to (M,\go)$ defined by
$$
\Phi_t (x) = \left( f_{2t} \circ \gf \right) (e^{-t} x),
\quad
\Psi_t (x) = \left( f_{2t} \circ \psi \right) (e^{-t} x) .
%\qquad t \ge 0, \; x \in B^{2n}(b) .
$$
Since $\left( f_t \right)^* \go = e^t \go$,
the embeddings $\Phi_t, \Psi_t$ are symplectic.
Moreover, $\Phi_0 = \gf$ and $\Psi_0 = \psi$.
For $T>0$ so large that $e^{-T} B^{2n}(b) \subset B^{2n}(b')$,
we have $\Phi_T = \Psi_T$.
Concatenating the path of embeddings $\Phi_t$, $t \in [0,T]$,
from $\gf$ to $\Phi_T$
with the path of embeddings $\Psi_{T-t}$, $t \in [0,T]$,
from $\Phi_T = \Psi_T$ to~$\psi$,
we obtain a required path of symplectic charts from
$\gf$ to~$\psi$.
\proofend

{\ni \it Remark.}
In the case where $(M^{2n},\go)$ is a cotangent bundle
$\left( T^*Q, d\gl\right)$,
a parametric version of the above argument gives a description of
the homotopy type of the space $\Emb \left( B^{2n}(b), T^*Q \right)$:
the map $\Emb \left( B^{2n}(b), T^*Q \right)\to Q$
defined by projecting the center of the ball to the base
is a Serre fibration with fibre homotopy equivalent to~$U(n)$.

\m
\proofof{Theorem~\ref{t:Liouville}}
If we prove that~$(M,\go)$ is symplectically aspherical,
then the theorem will follow from Proposition~\ref{p:emb} and
Theorem~\ref{t:rho}.
Let $X$ be a forward complete Liouville field on~$(M,\go)$.
Denote by $\{f_t\}$, $t \ge 0$, its forward flow.
Let $g \colon S^2\to M$ be a smooth map.
Denote $g_t=f_t\circ g$.
Since $\go$ is closed and all maps $g_t$ are homotopic, we have
\[
 \int_{S^2} g^*\go \, =
 \int_{S^2} g_t^*\go \,= \int_{S^2} g^* (f^*_t\go)\,= e^t \int_{S^2} g^*\go
\]
for each $t>0$. Thus $\int_{S^2} g^*\go$  vanishes,
and $(M,\go)$ is symplectically aspherical.
\proofend

%\subsection{Camel space}
If the space $\Emb \left( B^{2n}(b), M,\go \right)$ is not connected,
the classification of product tori
can be more complicated:

\begin{example}  \label{ex:camel}
{\rm
The {\it camel space}\/ with eye of size~$c>0$ is the open subset
$$
\cc^{2n}(c) \,=\, \{ x_1<0\} \cup \{ x_1>0\} \cup \Bcirc^{2n}(c)
$$
of $(\RR^{2n},\go_n)$.
Fix $b>0$ and define the symplectic embeddings
$\gf_{\pm} \colon B^{2n}(b) \to \cc^{2n}(c)$
by
$$
\gf_{\pm} (x_1,y_1, \dots,x_n,y_n) \,=\, \left( x_1 \pm \sqrt{b/\pi},y_1, \dots, x_n,y_n \right) .
$$
If $b \ge c$, then the maps $\gf_{\pm}$ are not
homotopic through symplectic embeddings,
see~\cite{EG,MT,V},
and hence $\Emb \left( B^{2n}(b), \cc^{2n}(c),\go_n \right)$ has at least $2$~components.
Let $\ba\in \RR_+^{2n}\?$ be such that $T(\ba) \subset B^{2n}(b)$.
The symplectomorphism
$$
\left( x_1,y_1,\dots,x_n,y_n \right) \,\mapsto\, \left( -x_1,-y_1, x_2,y_2,\dots,x_n,y_n \right)
$$
maps $\gf_- \bigl( T(\ba) \bigr)$ to
$\gf_+ \bigl( T(\ba) \bigr)$,
and hence
$\gf_- \bigl( T(\ba) \bigr) \sim \gf_+ \bigl( T(\ba) \bigr)$.
However, if $\ba$ is such that $\ua \ge c$, then
$\gf_- \bigl( T(\ba) \bigr) \not\approx \gf_+ \bigl( T(\ba) \bigr)$
by the Lagrangian Camel Theorem of~\cite{Th-99}.
Therefore, the connectedness requirement cannot be omitted in Theorem~\ref{t:rho}.
The classification of product tori in $\cc^{2n}(c)$ up to
Hamiltonian isotopy may be difficult.
Indeed, there might exist a symplectic embedding
$\gf \colon B^{2n}(b) \to \cc^{2n}(c)$ whose image is so tangled up
in the eye of $\cc^{2n}(c)$ that $\gf \bigl( T(\ba) \bigr)$
is Hamiltonian isotopic to neither of
$\gf_{\pm} \bigl( T(\ba) \bigr)$.
}
\end{example}

%%%%%%%%%%%%%%%%%%%%%%%%%%%%%%%%%%%%%%%%%%%%%%%%%%%%%%%%%%%%%%%%%%%%%%%%%%%%%%%%%%%%%%%%

\section{Proof of Theorem~\ref{t:shift}}  \label{s:shift}

\subsection{}

Consider symplectic polar coordinates $(\rho,\theta)$ 
on $\dot \RR^2 := \RR^2 \setminus \{0\}$ 
defined by
$$
(x,y) \,=\, \left( \sqrt{\rho/\pi} \cos 2\pi \theta, \sqrt{\rho/\pi} \sin 2\pi \theta \right),
\quad \rho>0,\,\, \theta \in S^1 = \RR/\ZZ.
$$
%.
%
For $s \in \RR$ and $m \in \ZZ$, define the domain 
\begin{equation*} \label{e:dms}
\cd_{m,s} \.=\. \left\{ \.(\rho_1,\theta_1,\rho_2,\theta_2) \mid \rho_2+s > m \rho_1\.  \right\}
\.\subset\.\RR^4 
\end{equation*}
and the map
$\Psi_{\?\?m,s} \colon \cd_{m,s} \to  \RR^4$,
\begin{equation*} \label{e:psims}
\Psi_{\?\?m,s} (\rho_1,\theta_1,\rho_2, \theta_2) \.=\. 
( \rho_1, \.\theta_1 + m \theta_2, \.\rho_2 + s - m \rho_1, \theta_2) . 
\end{equation*}
The map $\Psi_{\?\?m,s}$ is a smooth symplectic embedding
(for the same reasons as the map $\Psi$ in the proof of Lemma~\ref{l:step1}).

Let $(M,\go)$ be a symplectic manifold,
and let $\gf \colon B^{2n}(b_+) \to (M,\go)$ be a symplectic chart. 
We denote by~$0_{2j}$ the origin in~$\RR^{2j}$.
The key step in the proof of Theorem~\ref{t:shift} is the following proposition.

\begin{proposition} \label{p:shear}
Let $k\ge1$, $d_1,\ldots,d_k, b_+ >0$. 
Let $S \in \pi_2(M)$ be such that $s:=\gs(S)$ is positive and
$$
d_1+\dots +d_k +s < b_+.
$$
Then there exist a neighbourhood~$U_k$  of the isotropic $k$-torus 
 $$
T^k_{\mathtt{i}}(d_1,\dots,d_k):=0_{2n-2k-2}\times T(d_1,\ldots,d_{k-\?1})\times 0_2 \times T(d_k)
$$ 
in the open ball $\Bcirc^{2n}(b_+\?)$
and a Hamiltonian symplectomorphism~$\psi_k$ of $(M,\go)$ such that 
$(\psi_k \circ \gf)(U_k) \subset \Bcirc_\gf^{2n}(b_+\?)$ 
and  the map $\psi^\gf_k := \gf^{-1} \circ \psi_k \circ \gf$ coincides with 
$\id_{2n-4}\times \Psi_{\?\?m,s} $ on~$U_k$, where $m = c_1(S)$.
\end{proposition}
We will need the following lemma.
\begin{lemma} \label{l:squeeze}
Given  positive numbers $d_1, \dots,d_{k-\?1}$, 
for each $\eps>0$ 
there is a Hamiltonian flow $\{\Xi_t\}$, $t\in [0,1]$, on $\RR^{2k}$
such that $\.\.\Xi_1$  maps the torus 
$$
T=T(d_1,\dots,d_{k-1})\times  0_{2}
$$ 
into $\bigl(\Bcirc^{2}(\eps)\bigr)^k$
and $\.\.\Xi_t$ maps $T$ into 
$
\Bcirc^2(d_1+\eps)\times\dots \times \Bcirc^2(d_{k-\?1}+\eps)\times \Bcirc^{2}(\eps)$ 
for all~$t$.
\end{lemma}

\proof
We start with the following 
\begin{lemma} \label{l:squeeze'}
Given a  positive number $d>0$, 
for each $\eps_0>0$ 
there exist $\delta= \delta(d,\eps_0)>0$ and
a Hamiltonian flow $\{\Xi\.^{d,\eps_0}_t\}$, $t\in [0,1]$, on $\RR^{4}$
with the following properties:

$\Xi\.^{d,\eps_0}_t\?\?$  maps 
$ T(d)\times  \Bcirc^2(\delta) $
into
 $\Bcirc^{2}(d+ \eps_0)\times \Bcirc^{2}(\eps_0)$  for all~$t\in [0,1]$;

$\Xi\.^{d,\eps_0}_1\?\?$  maps
$T(d)\times  \Bcirc^2(\delta) $
into $\Bcirc^{2}(\eps_0)\times \Bcirc^{2}(\eps_0)$.
\end{lemma}
\proof
%F in this proof I have replaced m by \ell, since m has different meaning in the overall section
For each $t \in [0,1]$ 
%F added
and for $\ell \in \NN$, 
define the map $E_{t,\ell} \colon S^1 \to\CC^2 =\RR^4$ by 
$$
E_{t,\ell}(\theta) \,=\,\Bigl(\,\sqrt{(1-t)\.d/\pi} \,\, e^{2\pi i \. \theta},\,
\sqrt{t\.d/(\ell \pi)}\, \,e^{2\pi i \ell\.\theta}
\Bigr).
$$
Then $E_{0,\ell}$ is a diffeomorphism onto~$T(d)\times 0_2$.
For  $t <1$, the map $E_{t,\ell}$ is an embedding because its first component is.  
The integral over $S^1$ of the $1$-form $E_{t,\ell}^*\lambda$,
where $\lambda=x_1\.dy_1 + x_2\.dy_2$ is a primitive of~$\go_2$,
does not depend on~$t$ because
$$
\int_{S^1} E_{t,\ell}^*\lambda
\,=\, 
\int_{S^1} E_{t,\ell}^*(x_1\. dy_1) +   \int_{S^1} E_{t,\ell}^*(x_2\. dy_2) 
\,=\,  (1-t)\.d +  t\.d  \,=\, d.
$$
It follows that for each $q\in\left]0,1\right[$ there is a Hamiltonian flow
$\{\Phi^{q,\ell}_t\}$, $t\in [0,1]$, such that 
%F $\Phi^{q,\ell}_t (T)            =E_{q\. t,\ell}( T(d)\times 0_2)$ 
$\Phi^{q,\ell}_t (T(d) \times 0_2) =E_{q\. t,\ell}( T(d)\times 0_2)$ 
for all $t\in [0,1]$.
The absolute value of the first component of the map $E_{t,\ell}$ is
decreasing with respect to~$t$; 
the second component of $E_{t,\ell}$ tends uniformly  to zero  as $\ell \to\infty$.
Therefore,  after choosing $\ell$ large enough, we can assume that 
the tori $E_{t,\ell}(T(d)\times 0_2)$ are contained in 
$B^{2}(d)\times \Bcirc^{2}(\eps_0)$ for all $t\in [0,1]$
and that the torus  $E_{1,\ell}(T(d)\times 0_2)$ 
is contained in~$\Bcirc^{2}(\eps_0)\times \Bcirc^{2}(\eps_0)$.
Then, after choosing $q$ sufficiently close to~$1$, we can achieve that  the torus
$E_{{q},\ell}(T(d)\times 0_2) =\Phi^{q,\ell}_1 (T(d)\times 0_2)$ is contained in
$\Bcirc^{2}(\eps_0)\times \Bcirc^{2}(\eps_0)$ as well.
Let $\{\Xi^{\.d,\eps_0}_t=\Phi^{q,\ell}_t\}$.
By continuity, there exists $\delta= \delta (d,\eps_0)>0$ such that
 $\Xi\.^{d,\eps_0}_t$ maps $T(d)\times \Bcirc^{2}(\delta)$ into
 $\Bcirc^{2}(d+ \eps_0)\times \Bcirc^{2}(\eps_0)$  for all~$t\in [0,1]$, and
$\Xi\.^{d,\eps_0}_1\?\?$  maps
$T(d)\times  \Bcirc^2(\delta) $
into $\Bcirc^{2}(\eps_0)\times \Bcirc^{2}(\eps_0)$.
\proofend

\ni
If $k=2$, then Lemma~\ref{l:squeeze} immediately follows from Lemma~\ref{l:squeeze'}.
Otherwise, applying Lemma~\ref{l:squeeze'} $k-1$ times,
we construct positive numbers
$$
\eps_1=\min\bigl(\delta(d_{k-\?1}, \eps),\eps\bigr), 
\,\,\eps_2=\min\bigl(\delta(d_{k-\?2}, \eps_1),\eps\bigr), 
\,\dots,\,\, \eps_{k-\?1}=\min\bigl(\delta(d_1, \eps_{k-\?2}),\eps\bigr)
$$ 
and Hamiltonian flows  
$\{\Xi\.^{d_{k-\?1},\.\.\eps}_t\},\{\Xi\.^{d_{k-\?2},\.\.\eps_1}_t\},
\dots, \{\Xi\.^{d_1,\.\. \eps_{k-\?2}}_t\}$ with the prescribed properties.
Consider the Hamiltonian flows  $\{\Phi^1_t\}, \{\Phi^2_t\}, \dots,\{\Phi^{k-\?1}_t\}$ 
on $\RR^{2k}$ such that 
$$
\Phi^1_t=\id_{2k-4}\times \.\.\Xi\.^{d_{k-\?1},\.\.\eps}_t \?,\,\,
\Phi^2_t=\id_{2k-6}\times \.\.\Xi\.^{d_{k-\?2},\.\.\eps_1}_t\!\?\times \id_{2},\.\dots\.,\.
\Phi^{k-\?1}_t\?\?= \.\.\Xi\.^{d_1,\.\. \eps_{k-\?2}}_t\!\?\times \id_{2k-4}  .
$$
For each $j\in\{1, \dots,k-\?1\}$, we have
%F here and also 2 lines later, the brackets changed (the factor \Bcirc^2(\eps_j) should not be touched by \Phi^j_t) 
$$
\Phi^j_t \bigl( T(d_1,\dots,d_{k-j})\?\times \?\?\Bcirc^2(\eps_j) \bigr) \?\times \?\?\bigl(\Bcirc^2(\eps)\bigr)^j \subset
T(d_1,\dots, d_{k-j-\?1})\times \Bcirc^2(d_{k-j}+\eps_{j-\?1})\times \bigl(\Bcirc^2(\eps)\bigr)^j
$$
for all~$t\in [0,1]$, 
and 
$$
\Phi^j_1 \bigl(T(d_1,\dots,d_{k-j})\?\times \?\?\Bcirc^2(\eps_j) \bigr) \?\times \?\?\bigl(\Bcirc^2(\eps)\bigr)^j \subset
T(d_1,\dots, d_{k-j-\?1})\times \Bcirc^2(\eps_{j-\?1})\times \bigl(\Bcirc^2(\eps)\bigr)^j,
$$
where $\eps_0=\eps$. 
Concatenating the flows $\{\Phi^1_t\}, \{\Phi^2_t\}, \dots,\{\Phi^{k-\?1}_t\}$
(and reparametrizing the result  to make it smoothly depending on~$t$),
we obtain the required flow~ $\{\Xi_t\}$.
\proofend

\subsection{}{\it Proof of Proposition~\ref{p:shear} for $k=1$.\ }
Denote $\cd=\RR^{2n-4}\times \cd_{m,s}$, 
 $$
\Psi\,=\,\id_{2n-4}\times \Psi_{\?\?m,s}\colon \cd\to \RR^{2n}.
$$
Let $e_1=d_1+s$.
Consider the maps $f_0,f_1\colon S^1\to  \RR^{2n},$
$$
f_0(\zeta)= 0_{2n-2}\times \bigl(\rho=d_1, \.\theta=\zeta\bigr),
\,\,\,\,
f_1(\zeta)= 0_{2n-2}\times \bigl(\rho=e_1, \.\theta=\zeta\bigr).
$$
We have $T^1_{\mathtt{i}}(d_1)=f_0(S^1)$,  $T^1_{\mathtt{i}}(e_1)=f_1(S^1)$, 
and  $\Psi\circ f_0=f_1$.
Let $f^{\gf}_0=\gf\circ f_0$,  $f^{\gf}_1=\gf\circ f_1$.

First we prove that there is 
$\hat\psi_1\in \Ham(M,\go)$ such that $\hat\psi_1\circ f^{\gf}_0=f^{\gf}_1$. 
Denote $Z=[0,1] \times  S^1$.
Consider the map $F\colon Z\to \RR^{2n},$
$$
F(v,\zeta)= 0_{2n-2}\times \bigl(\rho=d_1+ v\.s, \.\theta=\zeta\bigr).
$$
We have  $f_0= F(0,\cdot\.)$, $f_1= F(1,\cdot\.)$, and
$$
\int_{Z} (\varphi\circ F)^*\omega
\, =
 \int_{Z} F^*\omega_n  
\, =
\int_{S_1} \!f_1^* (\rho\. d \theta) - \int_{S_1} \!f_0^* (\rho\. d \theta) \,=\, s.
$$
Taking the connected sum of $\varphi\circ F$ with a map $S^2\to M $ 
representing the class $-S$, we obtain a smooth map  
$\widehat F\colon Z\to M$ such that  $\widehat F$
coincides with $\varphi\circ F$  at the boundary of $Z$
(that is, $f^{\gf}_0=  \widehat F(0,\cdot\.)$, $f^{\gf}_1=  \widehat F(1,\cdot\.)$) and
$$
\int_Z   \widehat F^*\omega =0.
$$
Then, according to~\cite[Appendix A]{LO}, there exists a Hamiltonian flow
$\{\hat\psi_t\}$ on $(M,\go)$ such that  the map 
$$
 \widetilde F \colon Z\to M, \,\,\,\, (v,\zeta) \mapsto \hat\psi_v (f^{\gf}_0(\zeta))
$$
is homotopic to $\widehat F$ relative to the boundary. 
In particular, this implies
$$
\hat\psi_1\circ f^{\gf}_0=f^{\gf}_1=\gf\circ\Psi\circ f_0,
$$
as required.
It follows that 
$\gf^{-1}\?\?\circ\hat\psi_1\circ \gf|_{T^1_{\mathtt{i}}(d_1)} =\Psi|_{T^1_{\mathtt{i}}(d_1)}$.
Pick a neighbourhood $W\subset \Bcirc^{2n}(b_+\?)$ of the circle  $T^1_{\mathtt{i}}(d_1)$ 
 such that the maps  
$\psi_W:= \gf^{-1}\circ\hat\psi_1\circ \gf|_W^{\vphantom1} $  
and $\Psi |_W^{\vphantom1} $ are well defined. 
We shall prove that there is a Hamiltonian symplectomorphism $\Phi$ with support in $W$
  and a neighbourhood 
$U_1$ of the circle $T^1_{\mathtt{i}}(d_1)$ in $W$ such that
\begin{equation}  \label{e:Phi}
  \Phi|_{U_1}^{\vphantom1} \,=\,\psi_W^{-1}\?\?\circ \Psi|_{U_1}^{\vphantom1} .
\end{equation}
Then the symplectomorphism $\psi_1\in \Ham(M,\go)$ that coincides with 
$\hat\psi_1\circ \gf\circ\Phi\circ\gf^{-1}$ 
on $\gf(W)$ and coincides with $\hat\psi_1$ outside $\gf(W)$ will
satisfy  $\gf^{-1} \circ \psi_1 \circ \gf|_{U_1}^{\vphantom1} =\Psi|_{U_1}^{\vphantom1}$
as required. 

Trivialize the tangent bundle of $\RR^{2n-2}\times\dot \RR^2$  using  the symplectic frame 
$$
  \xi  =
\bigl(\partial _{x_1},\partial _{y_1}, \dots,
\partial _{x_{n-\?1}},\partial _{y_{n-\?1}},\partial _{\rho_n},\partial _{\theta_n}
\bigr).
$$
Denote by 
$$
\eta_w\colon \RR^{2n}\?\?\to T_w\bigl(\RR^{2n-2}\times\dot \RR^2\bigr),\,\,\,\, 
w\in \RR^{2n-2}\?\?\times\dot \RR^2
$$
the corresponding trivialization maps.
Let $\Sp(2n)$ denote the group of linear symplectomorphisms of~$\RR^{2n}$.
Consider the loop
$$
   g\colon S^1 \to \Sp(2n), \,\,\,\,
   g(\zeta) \,=\, 
 \eta_{f_0(\zeta)}^{-1}\circ d\.\bigl(\psi_W^{-1}\?\?\circ \Psi \Bigr)
   \circ \eta_{f_0(\zeta)}^{\vphantom1}.
$$
Recall that the fundamental group of $\Sp(2n)$ is isomorphic to~$\ZZ$; 
this gives rise to a function $\mu$ called the {\it Maslov  index\/} assigning 
to each continuous map $ S^1 \to \Sp(2n)$ an integer (see~\cite[p.\hskip.1em 48]{MS}).
\begin{lemma} \label{l:maslov}
The Maslov index of $g$ vanishes.
\end{lemma}
\proof
Define the maps $g_0,g_1\colon S^1 \to \Sp(2n)$,
$$
    g_0(\zeta) \,=\,  \eta_{f_1(\zeta)}^{-1}\circ d \Psi \circ \eta_{f_0(\zeta)}^{\vphantom1},
\,\,\,\,
    g_1(\zeta) \,=\,  \eta_{f_1(\zeta)}^{-1}\circ d \psi_W \circ \eta_{f_0(\zeta)}^{\vphantom1}.
$$
Since $\mu$ is additive with respect to the multiplication in $\Sp(2n)$~\cite[Theorem~2.29]{MS},
we have $\mu(g)= \mu(g_0)-\mu(g_1)$.
By the definition of~$\Psi$, we have $g_0(\zeta)=\id_{2n-4}\times A_\zeta \times \id_{2}$,
where $A_\zeta$ acts on $\CC=\RR^2$ as complex multiplication by $e^{2\pi i m \zeta}$.
Hence, according to \cite[p.\hskip.1em 49]{MS}, $\mu(g_0)=m$.

In order to compute the Maslov index of $g_1$, 
consider the torus $K$ constructed from  two copies, $\Sigma_1$ and $\Sigma_2$, of the annulus $Z=[0,1] \times  S^1$
by gluing together the respective boundary components.
Define the map $u\colon K\to M$ that coincides with $\varphi\circ F$ on $\Sigma_1$, 
and with $\widehat F$ on~$\Sigma_2$.
Orient $K$ by the volume form $dv\wedge d\zeta$ on~$\Sigma_2$.
Then the homology class of $u(K)$ is~$S$.
Consider the symplectic vector bundle  $u^*TM$ over~$K$.
Trivialize it over $\Sigma_1$ by means of the frame $\varphi_*\xi$, 
and over $\Sigma_2$, at the point $(v,\zeta)$, by means of the frame $(\hat\psi_v\circ\varphi)_*\xi$.
Then it follows from \cite[p.\hskip.1em 75]{MS} that $\mu(g_1)= c_1(u(K))=m$. 
Hence $\mu(g)=0$.
\proofend

Denote by $\Sp_1(2n)$ the subgroup of the group $\Sp(2n)$ consisting of the maps sending the 
vector $(0,\ldots,0,1)$ to itself.
The loop $g$ takes values in~$\Sp_1(2n)$.
By Lemma~\ref{l:maslov}, $g$ is contractible in~$\Sp(2n)$.
We claim that it is also contractible in~$\Sp_1(2n)$. 
Indeed, the inclusion $i\colon\Sp_1(2n)\hookrightarrow \Sp(2n)$ is the fiber of the smooth fibration
$$
  \pi\colon \Sp(2n) \to   \RR^{2n}\?\?\setminus\?\?\{0\},  \,\,\,\,  A \mapsto A (0,\ldots,0,1) .   
$$
It follows from the long exact sequence  of $\pi$ that $i$ induces an isomorphism of fundamental groups when $n\ge2$.
Thus there is a smooth family of maps $g^t\colon S^1 \to \Sp_1(2n)$, $t\in [0,1]$, such that 
%F for coherence, all the id symbols are now in typewriter
$g^0=\id$ and $g^1=g$. 
  
There is a linear isomorphism $I$ from the space of quadratic forms on 
$\RR^{2n}$ to the Lie algebra $\mathtt{sp}(2n)$ of the Lie group  $\Sp(2n)$
that assigns to a quadratic form $h$ the  Hamiltonian vector field generated by~$h$.
The quadratic forms that vanish on the line $\{(0,\ldots,0,{\cdot})\}$
are isomorphically mapped by $I$ to  the Lie algebra $\mathtt{sp}_1(2n)$ of~$\Sp_1(2n)$.
From the family $\{g^t\}$ we construct a smooth family of Hamiltonian functions $\{H_t\}$ with support in $W$
such that 
$$ 
\eta_w^{-1}\left(\.d^{\.\.2}\.(H_t)\right)= I^{-1} \?\?\left(\dot g^t(\theta_n)\right)
$$ 
%F \eta_w is a trivialisation. Here it is applied to the Hession d^2 H_t. Do you mean ``the Hession written in the 
% coordinates defined by \eta_w''? Is this really necessary, or may we omit \eta_w^{-1} here?
% 
for all 
$w=(x_1,y_1, \dots, x_{n-\?1},y_{n-\?1},\rho_n,\theta_n) \in T^1_{\mathtt{i}}(d_1)$, $t\in[0,1]$.
%F This point w is actually $w= (0_{2n-2}, d_1, \theta_n)$. Maybe it is better to write this more explicit vector ?
Then the time 1 flow $\Phi_+$ generated by  $\{H_t\}$ fixes each point $w\in T^1_{\mathtt{i}}(d_1)$ and
has the same differential as $\psi_W^{-1}\?\?\circ \Psi$ at~$w$.  

The symplectomorphism 
%F in the next line, \Phi_1^{-1} replaced by \Phi_+^{-1}
$\Upsilon:= \Phi_+^{-1}\?\?\circ\psi_W^{-1}\?\?\circ \Psi$ 
fixes $T^1_{\mathtt{i}}(d_1)$ pointwise and
satisfies $d\Upsilon(w)=\id$ for all  $w\in T^1_{\mathtt{i}}(d_1)$.
We shall prove that
there is a Hamiltonian symplectomorphism $\Phi_1$
with support in $W$ coinciding with $\Upsilon$ near~$T^1_{\mathtt{i}}(d_1)$. 
Then 
%F $\Phi= \Phi_1\circ\Phi^+\?$ 
   $\Phi= \Phi_+\circ\Phi_1\?$ 
is as required.

To construct $\Phi_1$, we use generating functions (cf.~\cite[Section~48]{Arnold}, \cite[Appendix~A.1]{HZ}).
Consider the  graph $\Gamma\subset\RR^{2n}\?\?\times\RR^{2n}$  of the map~$\Upsilon$. 
Denote by $T^{\times}\subset \Gamma$ the circle consisting of the points $(w,w)$, where $w\in T^1_{\mathtt{i}}(d_1)$.
Denote by $p=(p_1,\dots,p_n)$, $q=(q_1,\dots,q_n)$ the symplectic coordinates on the first copy of~$\RR^{2n}$,
and by $p'=(p'_1,\dots,p'_n), q=(q'_1,\dots,q'_n)$ those on the second copy.
By construction, $\Gamma$ is tangent  to the diagonal $\Delta\subset \RR^{2n}\?\?\times\RR^{2n}$ along~$T^{\times}$. 
Hence there is a tubular neighbourhood~$V$ of $T^{\times}$ in~$\Gamma$ such that
the map
$$
  \tau\colon  V \to \RR^{2n}, \quad (p,q,p',q') \mapsto (p',q)
$$ 
is a diffeomorphism onto a neigbourhood $U$ of  $T^1_{\mathtt{i}}(d_1)$ in~$W$.
Since $\Upsilon$ is symplectic, $V$ is Lagrangian with respect to the symplectic form 
$$
\Omega=-dp\wedge dq + dp'\?\wedge dq' = dq\wedge dp + dp'\?\wedge dq'.
$$
The $1$-forms $\alpha= -p\.dq  + p' dq'$, $\alpha'= q\.dp  + p' dq'$ satisfy 
$d\alpha=d\alpha'=\Omega$ and $\alpha=\alpha'-d(pq)$.
Thus the restrictions of  $\alpha$ and $\alpha'$ to~$V$ are closed.
They are exact because the restriction of~$\alpha$ to the diagonal~$\Delta$,
and hence to the circle $T^{\times}\subset V\cap \Delta$, vanishes.
%F  six lines later, f denotes another function (the difference F - p'q). I thus replaced here f by h (twice)
Let $h\colon V\to \RR$ be a primitive of $\alpha'$. 
Define $F\colon \tau(V) \to \RR$, $F= h\circ \tau^{-1}$.
Then $F$ is a generating function for~$V$, namely, $V$ is given by the equations 
$$
    q =  \frac{\partial F(p',q)}{\partial p'}, \quad
    p' =  \frac{\partial F(p',q)}{\partial q}.
$$
Note that $p'q$ is a generating function for~$\Delta$.

Since $\Gamma$ is tangent to~$\Delta$ along~$T^{\times}$,
the functions $F(p',q)$ and $p'q$ have the same respective first and second
partial derivatives at the points of the circle $T^1_{\mathtt{i}}(d_1)=\tau(T^{\times})$.
Thus the function 
%F  f(p,q) changed to f(p',q)
$f(p',q) := F(p',q)-p'q$ is $C^2$ small near $T^1_{\mathtt{i}}(d_1)$, and 
there exists a family of $C^\infty$ smooth functions $f_\delta\colon \RR^{2n}\to \RR$, 
defined for sufficiently small positive~$\delta$,  
such that the  function $f_\delta$ has support in the $\delta$-neighbourhood $W_\delta$ of~$T^1_{\mathtt{i}}(d_1)$,
coincides with $f$ on a smaller neighbourhood of~$T^1_{\mathtt{i}}(d_1)$, and
tends to zero in the $C^2$ topology as $\delta$ tends to zero.
(To explicitly construct such a family, we can proceed as follows. 
Pick a family of smooth compactly supported functions $\lambda_\delta \colon \left[ 0,\delta\right [ \to \left[ 0,\delta\right [$  
such that $\lambda_\delta$ is identity near~$0$ and its first and second derivatives are bounded uniformly over~$\delta$.
Given $x \in W_\delta$, denote by $x_0$ the point of~$T^1_{\mathtt{i}}(d_1)$ closest to~$x$ and 
draw the ray starting at $x_0$ and passing through~$x$.
Let $G_\delta\colon W_\delta\to W_\delta$ be the map
that sends $x$ to  the point $y$   such that $y$  lies on this ray and $\mathrm{dist} (y,x_0)
=\lambda_\delta\left( \mathrm{dist} (x,x_0)\right )$. 
Define $f_\delta$ to coincide with $f\circ G_\delta$ on~$W_\delta$.)

Denote by $L_\delta^t$ the Lagrangian submanifold 
in $\RR^{2n}\?\?\times\RR^{2n}$ defined
by the generating function $p'q + t f_\delta (p',q)$.
Picking $\delta$ sufficiently small, we can assume that each  
of the manifolds $L_{\delta}^t$ is sufficiently $C^1$ close to $\Delta$ 
and hence is a graph of a compactly supported symplectomorphism~$\Phi_t$.
The symplectomorphism $\Phi_1$ is Hamiltonian because
$\Phi_0=\id$ and $H^1(\RR^{2n})=0$.
Making $\delta$ smaller if necessary, we can assume that each $\Phi_t$ has support in~$W$.
Since $p'q + f_\delta (p',q)$ coincides with $F$ near~$T^1_{\mathtt{i}}(d_1)$,
the 
%F contactomorphisms 
symplectomorphisms
$\Phi_1$ and $\Upsilon$ also coincide near~$T^1_{\mathtt{i}}(d_1)$.
Thus $\Phi_1$ is as required, which concludes the proof of Proposition~\ref{p:shear} for $k=1$.
\proofend

\subsection{}{\it Proof of Proposition~\ref{p:shear} for $k>1$.\ }
Applying Proposition~\ref{p:shear} for $k=1$ to  the circle 
$T^1_{\mathtt{i}}(d_k)$, we obtain a neighbourhood $U_1$ of 
$T^1_{\mathtt{i}}(d_k)$ and a Hamiltonian symplectomorphism $\psi_1$ 
such that  $\psi_1^\varphi |_{U_1}^{\vphantom1}= \Psi |_{U_1}^{\vphantom1}$.
We shall construct a neighbourhood $U_k\subset \cd$ of
the torus $T^k_{\mathtt{i}}:=T^k_{\mathtt{i}}(d_1,\dots,d_k)$
and Hamiltonian symplectomorphisms $\Theta,\Theta_{\?\star}$
with support in $\Bcirc^{2n}(b_+\?)$ such that 
%
%\begin{equation}\label{e:uk}
\[
\Theta(U_k)\subset U_1,\quad \,\,
\Psi \circ \Theta|_{U_k}=\.\.\Theta_{\?\star} \circ \Psi|_{U_k}.
\]
%\end{equation}
%
Denote by  $\Theta^{\varphi}$ (resp.~$\Theta_{\?\star}^{\varphi}$)
the Hamiltonian symplectomorphism of $(M,\omega)$ 
that coincides with $\varphi \circ \Theta \circ  \varphi^{-1}$ 
(resp.~$\varphi \circ \Theta_{\?\star} \circ \varphi^{-1}$) on $B^{2n}_\varphi(b_+\?)$ 
and with the identity elsewhere. 
The symplectomorphism 
$\psi_k=(\Theta_{\?\star}^{\varphi})^{-1}\circ \psi_1 \circ  \Theta^{\varphi}$
will then have the required property since
\[
             \varphi^{-1}\circ \psi_k \circ \varphi|_{U_k}^{\vphantom1}
 \,=\, 
               \Theta_{\?\star}^{-1}\circ  \psi_1^\varphi \circ \Theta|_{U_k}^{\vphantom1}
 \,=\, 
               \Theta_{\?\star}^{-1}\circ  \Psi \circ \Theta|_{U_k}^{\vphantom1} 
 \,=\, 
            \Psi|_{U_k}^{\vphantom1}.
\]

It  remains to construct $\Theta$ and~$\Theta_{\?\star}$. 
Let $\eps>0$.
Applying Lemma~\ref{l:squeeze}, we obtain a Hamiltonian flow $\{\Xi_t\}$
on $\RR^{2n}$ such that $\Xi_1$  maps the torus $T=T(d_1,\dots,d_{k-1})\times  0_{2}$ 
into  $\bigl(\Bcirc^{2}(\eps)\bigr)^k$ and 
\begin{equation}\label{e:Xit}
\Xi_t(T)\subset \Bcirc^2(d_1+\eps)\times\dots \times \Bcirc^2(d_{k-\?1}+\eps)\times \Bcirc^{2}(\eps) \,\,\,
                                                                                   \mbox{for all} \, \,\,t\in [0,1].
\end{equation}
Consider the Hamiltonian flow  
$$
\left\{\.\mathrm{P\?\?}_{t}=\id_{2n-2k-2}\times  \.\.\.\Xi_t\?\?\times \id_{2}  \right\}, \,\, t\in [0,1],
$$
on~$\RR^{2n}$. 
Let $b'=b_+\?\?\? - s$.
Clearly, the torus $T^k_{\mathtt{i}}$ is contained in $\cd \cap \Bcirc^{2n}(b')$. 
We claim that by choosing $\eps$ sufficiently small
we can achieve that $\mathrm{P\?\?}_{t}$ maps $T^k_{\mathtt{i}}$  
into $\cd \cap \Bcirc^{2n}(b')$ for all~$t\in [0,1]$, 
and that $\mathrm{P\?\?}_{1}$ maps $T^k_{\mathtt{i}}$ into~$U_1$. 
Indeed, if $m\.\eps< d_k$, then the set 
$$
\Bcirc^{2}(\eps) \times T(d_k)= \left\{\rho_1<\eps,\.\. \rho_2=d_k\right\}
$$
is contained in~$\cd_{m,s}$.
It follows from~(\ref{e:Xit}) that
for all $t$ the torus $\mathrm{P\?\?}_{t}(T^k_{\mathtt{i}})$ is contained in 
$\RR^{2n-4}\times \Bcirc^{2}(\eps) \times T(d_k)$, and hence in~$\cd$. 
If $d_1+\dots+d_k+k\.\eps<b'$, then it follows from~(\ref{e:Xit}) that
$\mathrm{P\?\?}_{t}(T^k_{\mathtt{i}}) \subset \Bcirc^{2n}(b')$ for all~$t$.
Finally, for $\eps $ such that $0_{2n-2k-2}\times\bigl(\Bcirc^{2}(\eps)\bigr)^k \times T(d_k)$ 
is a subset of~$U_1$,
we have  $\mathrm{P\?\?}_{1}(T^k_{\mathtt{i}}) \subset U_1$.

It follows from the definition of the map $\Psi$ that 
$\Psi(\mathrm{P\?\?}_{t}(T^k_{\mathtt{i}}))$ is contained in 
$\Bcirc^{2n}(b_+\?)$ for all $t\in [0,1]$.
Therefore, there is an open set $W\subset \cd \cap \Bcirc^{2n}(b')$
that contains all the tori $\Psi(\mathrm{P\?\?}_{t}(T^k_{\mathtt{i}}))$
and satisfies  $\Psi(W) \subset \Bcirc^{2n}(b_+\?)$.
Then there exists a neighbourhood $U_k$ of the torus  $T^k_{\mathtt{i}}$
such that $\mathrm{P\?\?}_{t}(U_k)\subset W$ for all~$t$, 
and $\mathrm{P\?\?}_{1}(U_k)\subset U_1$.

Applying to $\{\mathrm{P\?\?}_t\}$  an appropriate cut-off,
we construct a Hamiltonian flow  $\{\mathrm{P}'_{\?\?t}\}$,  $t\in [0,1]$,
with support in $W$ such that 
$\mathrm{P}'_{\?\?t}\.|_{U_k}^{\vphantom1}=\mathrm{P\?\?}_t\.|_{U_k}^{\vphantom1}$
for all~$t$ and $\mathrm{P}'_{\?\?1}(U_k)\subset U_1$.
Define the Hamiltonian flow 
$\{\mathrm{P}^\star_{\?\?t}\}$,  $t\in [0,1]$,
with support in $\Psi(W)\subset \Bcirc^{2n}(b_+\?)$ 
by $\mathrm{P}^\star_{\?\?t}= \Psi\circ \mathrm{P}'_{\?\?1} \circ \Psi^{-1}.$
Then $\Theta= \mathrm{P}'_{\?\?1}$ and~$\Theta_{\?\star}=\mathrm{P}^\star_{\?\?1}$
are as required.
\proofend

\subsection{}{\it Proof of Theorem~\ref{t:shift}.}
It suffices to prove the theorem
under the additional assumption that $d_j=e_j$ for $j<k$.
Indeed, in view of Theorem~\ref{t:t1}~(i), the claim will then also hold for vectors 
that differ at only one component; after that the general case follows
by changing one component at a time.

We extend the symplectic chart~$\gf$ from $B^{2n}(b)$ to a larger ball $B^{2n}(b_+)$
with $b_+\?\?>b$, and keep the letter~$\gf$ for this extension.
For $\bd'=(d'_1,\dots,d'_k)$, we abbreviate 
$T_\varphi (a,\dots,a,a+d'_1,\dots,a+d'_k)$ to  $ T_{\varphi\?,a}(\bd')$.
Given $\tau\in \bigl[0,\min(c,b_+\?\?-b)\bigr[$, 
denote by $\cv_\tau$ the subset of $\RR^k$ formed by vectors $(d_1,\dots,d_k)$
such that $d_1+\dots +d_k \le b+\tau$ and $d_j\ge c-\tau$ for all $j\in\{1,\dots,k\}$.
Pick $\delta\in \bigl ]0,\min(c,b_+\?\?-b)\bigr[$.
Recall that $\gs_a(S)=\gs(S) - c_1(S) a$.
\begin{lemma} \label{l:shift}
Let $S \in \pi_2(M)$.
There exists $A_S>0$ 
such that for each $a\in\left]\.0,A_S\right]$
and for each pair of vectors 
$$
\bd=\left(d_1,\dots,d_{k-\?1},d_k\right), \,\,\,\,
\bd_S=\left(d_1,\dots,d_{k-\?1},d_k+\gs_a(S)\.\right)
$$ 
belonging to $\cv_\delta$, we have 
$T_{\varphi,a}(\bd)\approx T_{\varphi,a}(\bd_S)$.
\end{lemma}
\proof
Denote $s=\sigma(S)$, $m=c_1(S)$.
Assume first that $s\ge0$.
It  follows from Proposition~\ref{p:shear} 
and the definition of the map $\Psi_{\?\?m,s}$
that  for each $\bd\in\cv_\delta$ there exist
a neighbourhood~$U$   of the isotropic $k$-torus 
 $T^k_{\mathtt{i}}(\bd)$ in $\Bcirc^{2n}(b_+)$
and a map $\psi\in \Ham(M,\go)$ such that 
for every torus $T(a_1,\dots,a_n)$ contained in $U$
we have 
$$
\psi \bigl(   T_\varphi (a_1,\dots, a_{n-\?1},a_n) \bigr)\.=\.
T_\varphi(a_1,\dots, a_{n-\?1},a_n + s -  m\.a_{n-\?1}).
$$
Therefore, by Theorem~\ref{t:t1}~(i), for each $\bd\in\cv_\delta$ there are
a positive number $A_{S\?,\.\bd}$ and a neighbourhood
$W_{S\?,\.\bd}$ of $\bd$ in $\cv_\delta$
such that for each  $\bd'\in W_{S\?,\.\bd}$ and each  $a\in\left]\.0,A_{S\?,\.\bd}\right]$
we have $T_{\varphi,a}(\bd')\approx T_{\varphi,a}(\bd'_S)$.

Since $\cv_\delta$ is compact, there are $\bd^{(1)},\dots, \bd^{(l)}\in \cv_\delta$
such that the sets $W_{S\?,\.\bd^{(j)}}$ cover  $\cv_\delta$.
Let $A_S$ be the smallest of the numbers $A_{S\?,\.\bd^{(j)}}$.
Then $T_{\varphi,a}(\bd)\approx T_{\varphi,a}(\bd_S)$
 for each  $\bd\in \cv_\delta$ and each  $a\in\left]\.0,A_S\right]$.
In particular,   $T_{\varphi,a}(\bd)\approx T_{\varphi,a}(\bd_S)$ 
for each $a\in\left]\.0,A_S\right]$ when $\bd,\bd_S\in \cv_\delta$.
The latter statement is invariant under changing the sign of~$S$,
and therefore we can drop the assumption that $s\ge0$.
\proofend

\m
Assume first that  $(M,\go)$ is not special. 
Let $S_1, \dots, S_r$ be elements of $\pi_2(M)$ 
such that their classes form the basis of the free Abelian group 
$\pi_2(M)/ \bigl( \ker \gs \cap \ker c_1 \bigr)$.
We can assume that $r\ge1$,
otherwise there is nothing to prove. 
Consider the free Abelian group $\sigma\left(\pi_2(M)\right)$.
If it is trivial, then $r=1$.
If its rank is $1$, then $r=1$ (otherwise  $(M,\go)$ would be special).
If the rank of this group is greater than~$1$,
then $r\ge 2$ and we can choose  $S_1, \dots, S_r$ such that for all $j\in\{1,\dots,r\}$ 
the numbers  $ s_j=\sigma(S_j)$  satisfy the inequality~$|s_j|<\delta$.
Denote $m_j=c_1(S_j)$. 
Pick $A>0$ such that for all  $j\in\{1,\dots,r\}$   we have 
\begin{equation*} \label{e:A}
A \le A_{S_j},    \,\,\,\,   | s_j-m_j A | < \delta .
\end{equation*}
If $(M,\go)$ is special, we set $r=1$, $S_1=S_0$ (or $S_1=-S_0$), and $A=A_{S_1}$.

Let $a\in\left]\.0,A\right]$.
Let
$$
\bd=\left(d_1,\dots,d_{k-\?1},d_k\right), \,\,\,\,
\be=\left(d_1,\dots,d_{k-\?1},e_k\right)
$$ 
be vectors in $\cv_\delta$.
We assume that the difference $d_k-e_k$
is an element of $G_a=\gs_a(\pi_2(M))$ if  $(M,\go)$ is not special,
and an element of  $G_a(S_0)=\gs_a(\langle S_0\rangle)$) if $(M,\go)$ is special.
Hence there are $n_1,\dots,n_r\in \ZZ$ such that 
$$
e_k-d_k\,=\, \sum_{j=1}^r n_j\.\gs_a(S_j)\,=\, \sum_{j=1}^r n_j\.(s_j-m_j a).
$$
After changing the signs of $S_j$ if necessary, we can assume that all 
coefficients $n_j$ are nonnegative. 
We need to prove that  $T_{\varphi,a}(\bd)\approx T_{\varphi,a}(\be)$.

Let $u_1,\dots,u_N$ be a sequence of numbers such that for each  
 $j\in\{1,\dots,r\}$ exactly $n_j$ of them equal $s_j-m_j a$.
It gives rise to the sequence $q_0,q_1,\dots,q_N$,
where $q_0=d_k$, $q_l=d_k+\sum_{i=1}^l u_i$ for all $l\in\{1,\dots,N\}$
(and hence $q_N=e_k$).
Without loss of generality, we can assume that $d_k<e_k$.
If 
\begin{equation} \label{e:int}
    \ q_l\in[\.d_k-\delta, e_k +\delta\.] \,\,\, 
\mbox{for all}  \,\,\,   l\in\{1,\dots,N\}  ,
\end{equation}
then  each of the vectors $\bq_l=(d_1,\dots,d_{k-\?1},q_l)$
belongs to~$\cv_\delta$.
Since  $a\le A_{S_j}$ for all $j$,  it then follows from Lemma~\ref{l:shift} that
\begin{equation*} 
T_{\varphi,a}(\bd) =T_{\varphi,a}(\bq_0) \approx T_{\varphi,a}(\bq_1)\approx
\dots
 \approx  T_{\varphi,a}(\bq_{N-1}) \approx T_{\varphi,a}(\bq_N) = T_{\varphi,a}(\be).
\end{equation*}
It remains to show that the sequence $u_1,\dots,u_N$
can be chosen to satisfy~(\ref{e:int}).
For $r=1$, there is no choice involved in the construction of the sequence,
and all $q_l$ belong to $[\.d_k, e_k ]$.
Let $r>1$. 
Then $|s_j-m_j a|<\delta$ for all $j$ since $|s_j|<\delta$ and $|s_j-m_j A|<\delta$.
We choose the numbers $u_l$ in succession,
using the following rule: if $q_{l-1}>e_k$, then $u_l<0$, 
and if $q_{l-1}<d_k$, then $u_l>0$.
Then (\ref{e:int}) will hold true.
This completes the proof of Theorem~\ref{t:shift}
\proofend

%%%%%%%%%%%%%%%%%%%%%%%%%%%%%%%%%%%%%%%%%%%%%%%%%%%%%%%%%%%%%%%%%%%

\appendix

\section{Areas of holomorphic curves in a hyperannulus} \label{a:annulus}

\renewcommand{\theequation}{A\arabic{equation}}
\setcounter{equation}{0}

For $r>0$, denote by $B_r$ (resp.~$\Bcirc_r$) the closed (resp.~open)
ball of radius $r$ in the complex vector space $\CC^n$ centred at the origin.
Denote $B_0=\{0\}$.
In this appendix, we prove the following

\begin{theorem}  \label{t::a}
Let $r_+\?\?>r_-\?\?\ge 0$.
Let $\.V\?\?$ be a holomorphic curve (a 1-dimensional analytic subvariety)
in the hyperannulus $\Bcirc_{r_+} \?\?\setminus B_{r_-}$
such that the closure of $\.\.V\?\?$  intersects~$\partial B_{r_-}$.
Then the area of $\.\.V\?\?$ is at least~$\pi \.(r_+^2\?\?-r_-^2)$.

If the area equals~$\pi \.(r_+^2\?\?-r_-^2)$, then  $\.\.V\?\?$ is the intersection
of a complex line in $\CC^n$ with the hyperannulus.
\end{theorem}
In the particular case where $r_-\?\?=0$, Theorem~\ref{t::a} is
equivalent to the $1$-dimensional version of the Lelong theorem
that gives a lower bound for the areas of holomorphic curves in
a ball passing through the centre~\cite{LeLong,V-89}.

Let $z_1=x_1\?+i\.y_1, \ldots, z_n = x_n\?+i\.y_n$
be the standard coordinates on~$\CC^n$.
Consider the $1$-form
$\alpha_n=\sum_{j=1}^n (x_j \,d y_j - y_j \,d x_j)$ on~$\CC^n$.

\begin{lemma}  \label{l::a}
Let $\gamma \colon S^1 \to \CC^n$ be a $C^1$-smooth curve.
Then its length $\ell(\gamma)$ satisfies the inequality
$$
\ell^2(\gamma) \.\. \ge \.\. 2\pi \?\? \int_{S^1} \! \gamma^* \alpha_n.
$$
\end{lemma}

\proof
For $n=1$, this is the classical isoperimetric inequality, see~\cite{O}.
The general case reduces to the $1$-dimensional case as follows.
Write $\gamma=(\rho_1, \ldots,\rho_n)$, where
$\rho_1, \ldots,\rho_n$ are maps from $S^1\?=\RR/\ZZ$  to~$\CC$.
By the Minkowski inequality (see~\cite{HLP}, p.\hskip.1em 146), we have
\begin{equation}
\label{e:a-Mi}
\begin{split}
\ell^2(\gamma) \. =  \left(\int_{\?S^1}    \! |  \gdot(t) | \. dt \right)^{\!2} 
 = 
  \left(\int_{\?S^1} \! \textstyle\sqrt{\sum_{j=1}^n |\dot\rho_j(t)|^2\.} \.\. dt \right)^{\!2} 
  \\
 \? \ge   \.
 \sum_{j=1}^n \! \left(\int_{\?S^1} \! |  \dot\rho_j(t)| \. dt \right)^{\!2}
\? = \. \sum_{j=1}^n  \ell^2(\rho_j),
\end{split}
\end{equation}
where $|\cdot|$ denotes the length of the vector,
and $\ell(\rho_j)$ is the length of the curve~$\rho_j$.
Since the isoperimetric inequality holds for the curves~$\rho_j$,
we have
$$
  \sum_{j=1}^n  \ell^2(\rho_j) \.\. \ge  \.\.
 2\pi \?\?  \sum_{j=1}^n \int_{S^1} \! \rho_j^* \alpha_1
 \.\. = \.\. 2\pi \?\? \int_{S^1} \!\? \gamma^* \alpha_n.
$$
\proofend

\proofof{Theorem~\ref{t::a}}
Let $Q$ be the subset of $\left] r_-,r_+\right[$ formed
by those~$r$ for which $V$ is transverse to~$\partial B_r$.
%F next 3 lines added
For reasons of analyticity,
the complement of $Q$ has no accumulation point in $\left] r_-,r_+\right[$.
In particular, $Q$ is of full measure in $\left] r_-,r_+\right[$.
Consider the function $F \colon \left[ r_-,r_+\right] \to \RR$,
where $F(r)$ is the area of the curve $V \cap (\Bcirc_{r} \setminus B_{r_-})$.
Then $F(r_-)=0$ and $F$ is monotone non-decreasing.
We shall prove that for each $r\in Q$,
the derivative $F'(r)$ exists and is not less than~$2 \pi r$.
Since $F$ is monotone non-decreasing, its derivative $F'$ is measurable and
$\int_{r_-}^{r_+} F'\,dr \,\le\, F(r_+)-F(r_-)$,
see Theorem 7.21 in~\cite{WZ}.
(Actually, it is easy to check that $F$ is continuous and hence the inequality is an equality.)
Therefore,
$$
F(r_+) \,=\, F(r_+)-F(r_-) \,\ge\, \int_{r_-}^{r_+}F'\, dr \,\ge\,
\int_{r_-}^{r_+}  \?\? 2 \pi r \.\.dr  \,=\, \pi \.(r_+^2\?\?-r_-^2).
$$

It follows from the maximum principle that $V$
intersects each sphere~$S_r=\partial B_{r}$.
Let $r\in Q$.
The set $V \cap S_r$ is the union of immersed circles
$W_1,\dots,W_m$, where $m\ge1 $.
Parametrize these circles by immersions
$\gamma_1,\dots,\gamma_m$ of $S^1\?=\RR/\ZZ$ into $S_r$
(such that, for each~$k$, the image of $\gamma_k$ is~$W_k$, and
$\gamma_k$ is an embedding outside finitely many points).
Consider the angle $\psi_k(t)\in \left[0, \pi/2\right]$ between the vector
$i\.\gdot_k(t)$ and the tangent space~$T_{\?\gamma_k\?\?(t)}S_r$.
Since $i\.\gdot_k(t)$ is tangent to $V$ and orthogonal to $\gdot_k(t)$,
 the angle between the tangent spaces to $V$ and
$S_r$ at the point $\gamma_k(t)$ also equals~$\psi_k(t)$.
Since $V$ intersects $S_r$ transversely,
this angle is always positive.
Denote by $s_k(t,\eps)$ the oriented (having the same sign as~$\eps$)
distance  between the point $\gamma_k(t)$ and the sphere
$S_{r+\eps}$, measured along the direction of $i\.\gdot_k(t)$.
We have $s_k(t,\eps) = \eps /\sin (\psi_k(t)) + O(\eps^2)$, and thus
$$
F'(r) \.\. = \.\. \sum_{k=1}^{m} \.\.\. \int_{S^1}   \?\?
         | \gdot_k(t) |   \frac{ \partial s_k(t,\eps) } {\partial \. \eps}   \.\. dt
       \.\.  = \.\. \sum_{k=1}^{m} \.\.\. \int_{S^1}
         \frac{| \gdot_k(t) |} {\sin (\psi_k(t))} \.\. dt.
$$
We shall prove that
\begin{equation}
\label{e:a-m}
 \int_{S^1}  \frac{| \gdot_k(t) |} {\sin (\psi_k(t))} \.\. dt \.\.\ge \.\. 2 \pi r
\end{equation}
for each $k$, and hence $F'(r) \ge 2\pi r m \ge 2 \. \pi r$.

Consider the Euler vector field
$\zeta=\sum_{j=1}^n (x_j \partial_{x_j} + y_j \partial_{y_j} )$
and the vector field
$\xi=\sum_{j=1}^n (x_j \partial_{y_j} - y_j \partial_{x_j} )$
on~$\CC^n$.
Let $q\in S_r$.
Multiplication by $i$ in the tangent space $T\?\?_q \CC^n$
takes the hyperplane $T\?\?_q \. S_r$
to the kernel $\ker_q \?\alpha_n$ of the
$1$-form $\alpha_n=\sum_{j=1}^n (x_j d y_j - y_j  d x_j)$,
and the vector $\zeta_q$ to the vector~$\xi_q\in T\?\?_q \. S_r$.
Since $\zeta_q$ is orthogonal to $T\?\?_q \. S_r$,
the vector $\xi_q$ is orthogonal to~$\ker_q \?\alpha_n$.
Fix $k\in \{1,\ldots,m\}$.
We write $\gamma$ for $\gamma_k$ and $\psi$ for~$\psi_k$.
Because multiplication by $i$ is an isometry,
the angle between $\gdot(t)$ and $\ker_{\gamma(t)} \?\alpha_n$ equals~$\psi(t)$.
Let $ u(t)$ be the non-zero vector in $T_{\gamma(t)}S_r$
obtained by projecting  $\gdot(t)$ along the hyperplane $\ker_{\gamma(t)} \?\alpha_n$ 
onto the line containing~$\xi_{\gamma(t)}$.
We can assume that $ u(t)$  is a positive multiple
of $\xi_{\gamma(t)}$ for all~$t$
(after reversing the orientation of $S^1$ if necessary).
Since $\alpha_n(\xi_q) = r^2\?\/ = r\. |\.\xi_q|$ for each $q\in S_r$,
we have $\alpha_n(u(t))=r\. |\.u(t)|$.
Thus
\begin{equation}
\label{e:a-a}
 \int_{S^1} \! \gamma^* \alpha_n \.\.= \.
 \int_{S^1} \! \alpha_n(u(t))\. dt\.\.= \.
 {r} \! \int_{S^1} \! |\.u(t)|\. dt .
\end{equation}
We have $|u(t)|=\sin (\psi(t)) \, |\gdot(t)|$, and hence
\begin{equation}
\label{e:a-b}
\int_{S^1}  \frac{| \gdot(t) |} {\sin (\psi(t))} \.\. dt
\.\. =
\int_{S^1}   \frac{|\gdot(t)|^2} {|\.u(t)|} \.\. dt.
\end{equation}
By the Cauchy--Schwarz inequality,
\begin{equation}
\label{e:a-c}
\int_{S^1}   \frac{|\gdot(t)|^2} {|\.u(t)|} \.\. dt \,
         \int_{S^1}  \!  |\.u(t)| \. dt \.\. \ge
\left(  \int_{S^1}  \!   |\gdot(t) | \. dt \right)^{\!\lo2}
 =\.\. \ell^2(\gamma),
\end{equation}
where $\ell(\gamma)$ is the length of~$\gamma$.
Applying Lemma~\ref{l::a} to $\gamma$ and using~\eqref{e:a-a},
we obtain the inequality
\begin{equation}
\label{e:a-d}
  \ell^2(\gamma) \.\. \ge \.\. 2\pi \?\? \int_{S^1} \! \gamma^* \alpha_n
  \.\. = \.\. 2\pi r \?\? \int_{S^1} \! |\.u(t)|\. dt.
\end{equation}
Combining~\eqref{e:a-b},~\eqref{e:a-c}, and~\eqref{e:a-d}, we get~\eqref{e:a-m}.
This proves the first statement.

In order to prove the second statement, suppose that
$\.\.V\?$ has area~$\pi \.(r_+^2\?\?-r_-^2)$.
Then all inequalities in the above argument, 
including those in the proof of Lemma~\ref{l::a}, 
turn into equalities.
We may assume that $\.\.V\?\?$ is irreducible, 
by considering only one of the irreducible components of $\.\.V\?$
(note that each component has positive area).

Pick $r\in Q$, and let $\gamma=(\rho_1, \ldots,\rho_n)$ be an immersion 
parametrizing a component of $V \cap S_r$. 
We may assume, after applying to $\.\.V\?$ a unitary transformation,
that the complex line tangent to $\.\.V\?$ at the point $\gamma(0)$
is parallel to the first coordinate axis. 
Then $\dot \rho_2(0)=\dots=\dot \rho_n(0)=0$.
The Minkowski inequality~\eqref{e:a-Mi} is an equality,
which implies that the functions $\dot\rho_1, \ldots,\dot\rho_n$ are proportional.
Since there is a point where all of them except the first one vanish,
the functions $\dot\rho_2, \ldots,\dot\rho_n$ vanish everywhere.
Therefore, the image of $\gamma$ is contained in a complex line $E$ 
parallel to the first coordinate axis.
Consider the intersection $E'$ of $E$ with the hyperannulus.
Because $E'$ contains the image of~$\gamma$, for reasons of analyticity
$E'$ is contained in the curve $\.\.V\?\?$.
But $\.\.V\?$ is irreducible, hence $\.\.V\?=E'$.
This completes the proof. 
\proofend

%%%%%%%%%%%%%%%%%%%%%%%%%%%%%%%%%%%%%%%%%%%%%%%%%%%%%%%%%%%%%%%%%%%

\section{Existence of low admissible paths}  \label{a:path}

\renewcommand{\theequation}{B\arabic{equation}}
\setcounter{equation}{0}

Let $k \ge 2$.
Given an ordered pair of different numbers $i,j\in \{1,\ldots,k\}$,
consider the operator $P_{ij}$
(resp.~$M_{ij}$, resp.~$I_{ij}$) in $\GL (k;\ZZ)$
that adds to the $i$-th component of a vector in $ \RR^k$ its $j$-th component
(resp.~subtracts from the $i$-th component the $j$-th component,
resp.~swaps the $i$-th component and the $j$-th component),
and does not change the other components.

Denote by $\RR_+$ the set of positive real numbers.
A sequence $\bd = \bd^0, \bd^1, \bd^2, \ldots, \bd^{\ell} = \be$
of vectors in $\RR_+^k$ is called an
{\it admissible path}\/ from $\bd$ to $\be$ if for each $s$
the vector $\bd^{s + 1}$ is obtained from $\bd^s$
by the action of one of the operators $P_{ij}, M_{ij}, I_{ij}$.
Given vectors $\bv, \bw \in \RR^k_+$, we write $\bv \le \bw$ if
there is a permutation $\gs$ of $\{1, \ldots, k\}$
such that $v_i \le w_{\gs (i)}$ for all $i \in\{ 1, \ldots, k\}$.
This defines a partial order on $\RR_+^k$.
We say that a path
$\bd = \bd^0, \bd^1, \ldots, \bd^\ell = \bd'$
is {\it low\/}
if for each $s \in\{\. 0, 1, \ldots, l\}$
we have $\bd^s \le \bd$ or $\bd^s \le \bd'$.
Given $\bu = (u_1, \ldots, u_k)\in \RR^k$, we write
\[
\langle \bu \.\rangle \. = \. \langle u_1, \ldots, u_k \rangle
\]
for the free Abelian subgroup in $\RR$ generated over $\ZZ$ by the numbers
$u_1, \ldots, u_k$.

The following theorem may be well known to specialists in number theory
or geometric group theory, but we were unable to find it in the literature.

\begin{theorem}  \label{t:low}
Given $\bd=(d_1,\ldots,d_k)$ and $\be=(e_1,\ldots,e_k)$ in $\RR_+^k$
such that
$\langle \bd \.\rangle = \langle \be \rangle$,
there is  a low admissible path from $\bd$ to~$\be$.
\end{theorem}

This appendix is devoted to the proof of Theorem~\ref{t:low}.
We start with two remarks concerning low admissible paths.
First, if the path
$\bd^0, \bd^1, \ldots, \bd^{\ell-1}, \bd^{\ell}$
is admissible, then the path
$\bd^{\ell}, \bd^{\ell-1}, \ldots, \bd^1, \bd^0$
is also admissible, because $M_{ij}^{-1} = P_{ij}$.
Second, the concatenation of a low path from
$\bd$ to $\bd'$ and a low path from $\bd'$ to $\bd''$
does not have to be low.
However, the concatenation is always low when $\bd'\le\bd$ or $\bd'\le\bd''$.

\begin{lemma}  \label{l:aux}
If $\bu, \bu' \in \RR^l$
and $\langle \bu \.\rangle = \langle \bu' \rangle$,
then there is $A\in \GL(l,\ZZ)$ such that $A \bu = \bu'\?$.
\end{lemma}

\proof
Consider the homomorphisms
$ h,h'\colon \ZZ^l \to \langle \bu\rangle,$
$ h (\bn)=(\bu;\bn),$ $ h' (\bn)=(\bu'\?;\bn),$
where $(\.\cdot\.\.;\cdot\.\.)$ is
the scalar product on~$\RR^l\?\?$.
Since $h$ and $h'$ are surjective homomorphisms of
free Abelian groups, there are splittings
$\ZZ^l\?\?= \ker (h) \oplus\Lambda = \ker (h')\oplus \Lambda'$,
where $\Lambda, \Lambda'$ are  subgroups of $ \ZZ^l\?\?$,
and the restrictions  of $h,h'$ to $\Lambda, \Lambda'$
respectively are isomorphisms onto~$\langle \bu\rangle$.
%F next three lines changed to get directions of the mappings right
Consider a homomorphism $B\colon \ZZ^l \to \ZZ^l\?\?$
such that $B{\mid}_{\Lambda'}=(h{\mid}_{\Lambda})^{-1}\?\circ h'{\mid}_{\Lambda'}$
and $B{\mid}_{\ker (h')}$ is an isomorphism from $\ker (h')$ onto~$\ker (h)$.
We have $B\in\GL(l,\ZZ)$ and $h'=h \circ B$.
Define $A\in\GL(l,\ZZ)$ to be the transpose of~$B$.
Then 
$$
(\bu'\?;\bn) \,=\, h'(\bn) \,=\, h(B\bn) \,=\, (\bu;B\.\bn)=(A\bu;\bn)
$$
for each~$\bn \in \ZZ^l$,
and hence $A \.\bu = \bu'\?$.
%
%Consider a homomorphism $B\colon \ZZ^l \to \ZZ^l\?\?$
%such that $B{\mid}_{\Lambda}=(h{\mid}_{\Lambda})^{-1}\?\circ h'{\mid}_{\Lambda'}$
%and  $B{\mid}_{\ker (h)}$ is an isomorphism from $\ker (h)$ onto~$\ker (h')$.
\proofend
\begin{lemma}  \label{l:rk1}
Given $\bd, \be \in \RR_+^k$ such that
$\langle \bd \.\rangle = \langle \be \rangle$
and $\langle\bd\.\rangle$ has rank~$1$,
there is a low admissible path from $\bd$ to $\be$.
\end{lemma}
\proof
Let $d_0>0$ be such that $\langle d_0 \rangle= \langle \bd \.\rangle$,
and let $\bd'=(d_0,\ldots,d_0)\in \RR_+^k$.
By repeatedly applying to $\bd$ and $\be$ operations~$M_{ij}$,
we construct admissible paths
$\bd,\bd^1,\ldots,\bd^\ell=\bd'$ and
$\be,\be^1,\ldots,\be^m=\bd'$,
where $\bd \ge \bd^1 \ge \dots \ge \bd^\ell$ and
$\be \ge \be^1 \ge \dots \ge \be^m$.
Then $\bd,\ldots,\bd',\ldots,\be$ is a low admissible path.
\proofend

\begin{lemma}  \label{l:rk2}
Given $\bd, \be \in \RR_+^2$ such that
$\langle \bd \.\rangle = \langle \be \rangle$
and  $\langle\bd\.\rangle$ has rank~$2$,
there is a low admissible path from $\bd$ to $\be$.
\end{lemma}

\proof
By Lemma~\ref{l:aux}, there exists $A\in\GL(2;\ZZ)$
such that $A(\bd)= \be$.
The following lemma then shows
that there is an admissible path from $\bd$ to~$\be$.

\begin{lemma}  \label{l:A2}
Given $\bd,\be \in \RR_+^2$ and $A \in \GL (2;\ZZ)$
such that $A (\bd) = \be$,
there is an admissible path from $\bd$ to~$\be$.
\end{lemma}

\proof
Let $\ce$ be the set of matrices in $\GL (2;\ZZ)$ with one entry~$0$
and three entries in~$\{1,-1\}$.
The set $\ce$ generates the group $\GL (2;\ZZ)$.
Indeed, $\GL (2;\ZZ)$ is generated by the three matrices
$$
P \,=\,
\begin{pmatrix}
1 & 1 \\
0 & 1
\end{pmatrix},
\qquad
I \,=\,
\begin{pmatrix}
0 & 1 \\
1 & 0
\end{pmatrix},
\qquad
Q_1 \,=\,
\begin{pmatrix}
-1 & 0 \\
0  & 1
\end{pmatrix}
$$
(see e.g.~\cite[p.\hskip.1em 43]{Harpe}),
and we have
$P \in \ce$,
$$
I = \begin{pmatrix}
0 &  1 \\
1 & -1
\end{pmatrix} P, \qquad
Q_1 = \begin{pmatrix}
-1 & 1 \\
0  & 1
\end{pmatrix} P.
$$
Moreover, $\ce$ is closed under taking inverses.
Hence we can write $A = E_\ell \cdots E_2 E_1$, where $E_s \in \ce$
for all $s \in \left\{ 1, \dots, \ell \right\}$.
Let $\bd^0 = \bd$, and let
$$
\widetilde\bd^s  \.=\.(\.x^s\?, y^s) \.\.= \. E_s \.{\cdots}\.\. E_2 E_1 (\bd),
\,\,\,\,\,\, \bd^s \.=\.(\.|\.\.x^s|, |\.\.y^s|\.)
$$
for each~$s>0$.
The  components of $\bd^s$ cannot vanish because they generate
the free Abelian group $\langle\bd\.\rangle$ of rank~$2$.
Hence $\bd^s\in \RR^2_+$ for each~$s$.
It suffices to construct for each
$s \in \left\{ 1, \dots, \ell \right\}$
an admissible path from $\bd^{\.s-\?1}$ to $\bd^{s}$.
Denote by $\ck$ the subset of $\GL (2;\ZZ)$ consisting of
diagonal matrices.
For each $s \in \left\{ 1, \dots, \ell \right\}$
there is $K_s \in \ck$ such that
$K_s (\widetilde\bd^s) = \bd^s$. Set $K_0= \id$.
Fix $s \in \{1, \dots, \ell\}$.
Denote $E_{s}' = K_{s} E_{s} K_{s-\?1}$.
Since $K_{s-\?1}^{-1} = K_{s-\?1}$,
we have $\bd^s =  E_{s}' (\bd^{\.s-\?1})$.
Since $\ce$ is invariant under
multiplication from the left and the right by elements of~$\ck$,
we have $E_{s}'\in\ce$.
Because $\bd^{\.s-\?1}\?, \bd^s \in \RR_2^+$,
each transformation $E\in \ce$ that maps
$\bd^{\.s-\?1}$ to $\bd^s$ coincides
(possibly after precomposing or/and postcomposing
it with the involution $(x,y) \stackrel{I}{\longmapsto} (y,x)$)
with one of following three:
\[
(b,c) \longmapsto (b+\?c,c), \quad
(b,c) \longmapsto (b\?-\?c,c), \quad
(b,c) \longmapsto (c\?-\?b,c).
\]
In particular, this is true for $E=E_{s}'$.
The path from $\bd^{\.s-\?1}$ to  $\bd^s$
can therefore be chosen to be one of the following admissible paths
(with a possible addition of the involution $I$ at the beginning
or/and at the end):
\begin{eqnarray*}
(b,c) &\stackrel{P_{12}}{\longmapsto}& (b+\?c,c), \\
(b,c) &\stackrel{M_{12}}{\longmapsto}& (b\?-\?c,c), \\
(b,c) &\stackrel{M_{21}}{\longmapsto}& (b,c\?-\?b)
     \,\stackrel{I}{\longmapsto}\, (c\?-\?b,b)
     \,\stackrel{P_{21}}{\longmapsto}\, (c\?-\?b,c).
\end{eqnarray*}
This proves Lemma~\ref{l:A2}.
\proofend

\m
We shall prove now that there exists a
low admissible path from $\bd$ to~$\be$.
We call an admissible path $\bd=\bd^0,\ldots,\bd^\ell=\be$
special if each of the moves from $\bd^{s-\?1}$ to $\bd^s$
is by one of the following three operators:
$P=P_{12}, M=M_{12}, I=I_{12}$.
Since $P_{21} = I  P_{12} I$ and $M_{21} = I  M_{12}  I$,
%it follows that 
every admissible path from $\bd$ to~$\be$
can be transformed into a special one.
In particular,  special admissible paths  exist.
Let $\frak p$ be a special admissible path  of minimal length.
We claim that $\frak p$ is low.

Assume first that some $P$ move precedes some
$M$ move in the path~$\frak p$.
Then the path must contain either the sequences of moves $P,M$
or the sequences of moves $P,I,M$
(note that adjacent $I$ moves cannot occur due to the minimality
of~$\frak p$).
In the first case, the sequence
\[
(c,d) \stackrel{P}{\longmapsto}
(c+d,d) \stackrel{M}{\longmapsto}  (c,d)
\]
can be removed from the path, in contradiction with the minimality
of~$\frak p$.
The second case,
\[
(c,d) \stackrel{P}{\longmapsto}  (c+d ,d) \stackrel{I}{\longmapsto}
(d,c+d) \stackrel{M}{\longmapsto}  (-c,d),
\]
is impossible because $-c$ is negative.
It remains to consider the case where for some $s$
each of the first $s$ moves  is either $M$ or $I$, and each
of the last $\ell-s$ moves is either $P$ or~$I$.
Then we have
\[
\bd \ge \bd^1 \ge \cdots \ge \bd^s \le \cdots \le \bd^{\ell-1} \le \be,
\]
and the path $\frak p$ is low.
\proofend

\begin{lemma}  \label{l:low-i}
Let $\bd=(d_1,\ldots,d_k)$, $\be=(e_1,\ldots,e_k)$ be
vectors in~$ \RR_+^k$ satisfying $\langle \bd \.\rangle =\langle \be \rangle$.
Assume that there is $i \in \{\.1,\dots,k\}$ such that
$d_i=e_i$, $d_i\le d_j$ for all~$j$, 
and $d_i$ is primitive (indivisible) in~$\langle \bd \.\rangle$.
Then there is a low admissible path from $\bd$ to~$\be$.
\end{lemma}

\proof
Assume for notational convenience that $i=k$.
By repeatedly subtracting the number $d_k =e_k$ from the
components of $\bd$ that exceed~$d_k$,
we construct an admissible path $\bd,\bd^1,\ldots,\bd^\ell$, where
$\bd \ge \bd^1 \ge \cdots \ge \bd^\ell$ and $\bd^\ell$ is such that
$d^\ell_k \ge d^\ell_j$ for all~$j$.
Using the same procedure, we obtain an admissible path
$\be, \be^1, \ldots, \be^m$, where
$\be \ge \be^1 \ge \cdots \ge \be^m$ and
$\be^m$ is such that $e^m_k \ge e^m_j$ for all~$j$.

Denote by $\pi$ the projection $\RR^k \to \RR^{k-\?1}$
that forgets the $k$-th component,
and define an equivalence relation on $\RR^k$ as follows.
Let $\bu=(u_1,\ldots,u_k)$, $\bv=(v_1,\ldots,v_k)$.
Then $\bu \sim \bv$ if
\begin{itemize}
\item[$\diamond$]
  $\langle \bu \.\rangle =\langle \bv \rangle$,
\item[$\diamond$]
  $\.u_k = v_k$,
\item[$\diamond$]
  there is an operator $A\in \GL(k-\?1,\ZZ)$
  such that all components of the vector $A(\pi(\bu))-\pi(\bv)$
  are integer multiples of~$u_k$.
\end{itemize}
We claim that $\bd \sim \be$.
Indeed, since $d_k$ is indivisible in~$\langle \bd \.\rangle$,
there is a subgroup
$\Lambda \subset \langle \bd \rangle$ such that
$\langle \bd \rangle = \langle d_k \rangle \oplus \Lambda$.
For each $j \in \{1,\dots,k-\?1\}$ there exists $d'_j \in \Lambda$
(resp.~$e'_j\in \Lambda$) that differs from $d_j$ (resp.~$e_j$)
by an integer multiple of~$d_k$.
Consider the vectors $\bd' = (d'_1,\ldots,d'_{k-\?1})$,
$\be' =(e'_1,\ldots,e'_{k-\?1})$
in~$\RR^{k-1}\?\?$.
The groups $\langle \bd' \rangle$ and $\langle \be' \rangle$
are subgroups of~$\Lambda$.
Since
$\langle d_k \rangle \oplus \Lambda
 = \langle d_k \rangle + \langle \bd' \rangle
 = \langle d_k \rangle + \langle \be' \rangle$,
we have $\langle \bd' \rangle = \langle \be' \rangle = \Lambda$.
By Lemma~\ref{l:aux}, there is $A \in \GL(k-1,\ZZ)$
that takes $\bd'$ to~$\be'$.
This $A$ has the required property.

By construction, we have
$\bd \sim \bd^\ell$ and
$\be \sim \be^m$.
Thus~$\bd^\ell {\sim} \be^m$.
For a positive number~$C$,
we say that an admissible path has height~${\le}\.\.C$ if
the components of every vector
involved in this path do not exceed~$C$.
We shall construct an admissible path of 
height ${\le}\.\.d_k$ from $\bd^\ell$ to~$\be^m$.
By concatenating the path $\bd, \ldots, \bd^\ell$,
the path from $\bd^\ell$ to~$\be^m$, and 
the path $\be^m, \ldots, \be$,
we then obtain a low admissible path from $\bd$ to~$\be$.

Let $A \in \GL(k-\?1,\ZZ)$ be such that the components
of the vector $A(\pi(\bd^\ell))-\pi(\be^m)$
are integer multiples of~$d_k$.
Denote by $Q_j$, $j \in\{\.1,\ldots,k-\?1\}$, the operator in
$\GL(k-1,\ZZ)$ that changes the sign of the~$j$-th component
and keeps all other components intact.
Each operator in $\GL(k-\?1,\ZZ)$ can be written as a product
of operators $Q_j$, $I_{ij}$, and~$P_{ij}$,
see e.g.~\cite[p.\hskip.1em 43]{Harpe}.
In particular, we have $A = A_r \cdots A_2 A_1$,
where either $A_s = Q_{j_s}$, or $A_s= I_{i_s j_s}$,
or $A_s = P_{i_s j_s}$ for each $s \in \{1,\dots, r\}$.
Denote $\bv^0 = \pi(\bd^\ell)$, $\bv^s = A_s \cdots A_1 (\bv^0)$.
Consider the map $\psi \colon \RR^{k-1} \to \RR_+^{k-1}$
that takes a vector $(x_1,\ldots,x_{k-\?1})$ to the vector
$(y_1,\ldots,y_{k-\?1})$,
where $y_j \in \left]0, d_k \right]$ and $x_j\?-y_j$
is an integer multiple of~$d_k$ for each~$j$.
Denote  $\bu^s=\psi\.( \bv^{s})$.
Then $\bu^0=\pi(\bd^\ell)$, $\bu^r=\pi(\be^m)$,
and, since $\psi \circ A_s\circ\psi = \psi \circ A_s$,
we have $\bu^s=\psi\.( A_s \bu^{s-\?1})$
for each~$s\in\{\.1,\dots,r \}$.

Denote by $\bw^s = (w_1^s,\ldots, w_{k-\?1}^s,d_k)$
the vector in $\RR_+^k$ obtained by complementing the
$k-\?1$ components of  the vector~$\bu^s$
with the $k$-th component equal to~$d_k$.
We have  $\bw^0=\bd^\ell$, $\bw^r=\be^m$, and $\bw^s\le d_k$ for each~$s$.
It suffices to prove that  there is an admissible path
of height ${\le}\.\. d_k$ from $ \bw^{s-\?1}$
to~$\bw^s$ for each~$s$.

If $A_s= I_{i j}$, then  $\bw^s= I_{i j}( \bw^{s-\?1})$ and
such a path consists of the single step~$I_{i j}$.
Let $A_s= Q_{j}$.
For notational convenience, assume that $k=2$, $j=1$.
Let $\bw^{s-\?1}=(c,d)$.
If $c=d$, then $\bw^{s}=\bw^{s-\?1}\?\?$,
and there is nothing to prove.
If $c<d$, then
$ \bw^{s}=(d\?-\?c,d)$, and the path
\begin{equation}
(c,d)\, \stackrel{M_{21}}{\longmapsto}\, (c,d\?-\?c)
     \, \stackrel{I_{12}}{\longmapsto}\, (d\?-\?c,c)
     \, \stackrel{P_{21}}{\longmapsto}\, (d\?-\?c,d)
\label{eq:d-c}
\end{equation}
is of height~$\le d$.

Let $A_s= P_{ij}$.
If $w^s_i+w^s_j \le d_k$, then  $\bw^s= P_{i j}( \bw^{s-\?1})$
and the required path consists of the step~$P_{i j}$.
Let $w^s_i+w^s_j > d_k$.
For notational convenience, assume that $k=3$, $i=1$, $j=2$.
We have $\bw^{s-\?1}=(b,c,d)$, where $b+c>d$.
If $c=d$, then  $\bw^{s}\?=\bw^{s-\?1}\?\?$,
and there is nothing to prove.
If $c<d$, then $ \bw^{s}\?=(b\?+\?c\?-\?d,c,d)$.
By~\eqref{eq:d-c}, there is an admissible path of height ${\le}\.\. d$
from $\bw^{s-\?1}$
to $\bw' = (b,d-c,d)$, and hence to $\bw'' = M_{12}(\bw) = (b+c\?-\?d, d\?-\?c,d)$.
By~\eqref{eq:d-c} again, there is an admissible path
of height ${\le}\.\. d$ from $\bw''$ to~$\bw^s$.
Lemma~\ref{l:low-i} is proved.
\proofend

\m
\proofof{Theorem~\ref{t:low}}
The proof is by induction on~$k$.
Lemma~\ref{l:rk1} and Lemma~\ref{l:rk2} prove the statement for $k=1$ and $k=2$.
We shall prove the statement for $k\ge3$ assuming that it holds for $k-\?1$.
In view of Lemma~\ref{l:rk1}, we can assume that $\rk \langle\bd\.\rangle\ge 2$.

\begin{lemma}  \label{l:indi}
Let $\bu=(u_1,\ldots,u_k)\in \RR_+^k$.
There is a low admissible path from  $\bu$ to a vector 
$\bu^+\?=(u^+_1,\ldots,u^+_k) \in \RR_+^k$ 
such that $\bu^+\?\le \bu$,  $u_k^+\?\le u_j^+$ for all~$j$, 
and $u^+_k$ is indivisible in  $ \langle\bu\rangle$.
\end{lemma}
\proof 
By repeatedly subtracting from the last component of the vector $\bu$
its other components, we construct a low admissible path from~$\bu$
to a vector $\bu'=(u_1,\dots,u_{k-\?1},u_k')\in \RR_+^k$ such that $\bu'\le \bu$
and $u_k'\?\le u_j$ for all~$j$.
If $u_k'$ is indivisible in  $ \langle\bu\rangle$, 
then  $\bu^+\?\?=\bu'$, and the lemma is proved.
Denote $\bu^-=(u_1,\dots,u_{k-\?1})$. 
If $u_k'$ is not primitive, then the rank~$m$ 
of the free Abelian group $\langle \bu^-\rangle$ is at least~$2$.
We claim that there is a basis $x_1,\dots,x_m$ of $\langle \bu^-\rangle$
formed by positive numbers smaller than~$u_k'$. 
Indeed, since  $\rk \langle \bu^-\rangle = m \ge 2$, there is a primitive
element $x_1\in \langle \bu^-\rangle$ satisfying $0< x_1<u_k'$.
Extend $x_1$ to a basis $x_1,\dots,x_m$ of $\langle \bu^-\rangle$,
and add to each $x_j$, $2\le j\le m$, a suitable integer multiple 
of $x_1$ to ensure that $0< x_j<x_1$.
After reordering the elements of the basis, we can assume that 
$x_m$ is the smallest of them.
Consider the vector $\bv^-=(x_1,\dots,x_1,x_2,\dots,x_m)$ in $\RR_+^{k-\?1},$
and the vectors
$\bv=(x_1,\dots,x_1,x_2,\dots,x_m,u_k')$, 
$\bu^+=(x_1,\dots,x_1,x_2,\dots,x_{m-1},u_k',x_m)$
 in $\RR_+^{k}$.
By the induction hypothesis, there is a low admissible path 
from $\bu^-$ to~$\bv^-$. 
Attaching to the vectors of this path $u_k'$ as the $k$-th component,
we construct a low admissible path from $\bu'$ to~$\bv$.
We have $\bv\le\bu' \le \bu$ and $\bu^+ \le \bu$. 
Therefore, concatenating the path from $\bu$ to~$\bu'$, 
the path from $\bu'$ to~$\bv$, and the one-step transposition path
from $\bv$ to~$\bu^+$, we obtain the required low admissible path. 
\proofend

By Lemma~\ref{l:indi}, we can assume that $d_k$ and $e_k$ 
are indivisible in~$\langle\bd\.\rangle$,
$d_k \le d_j$ and $e_k\le e_j$ for all~$j$. 
If $d_k=e_k$, then Lemma~\ref{l:low-i} (with $i=k$)
proves the statement.

Otherwise, $\langle d_k,e_k\rangle$ is a free Abelian group
of rank~$2$.
Consider the elements $x\in \langle\bd\.\rangle$ such that
$n_x x\in\langle d_k,e_k\rangle$ for some non-zero integer~$n_x$.
They form a rank~$2$ free Abelian subgroup~$\Delta$ of~$\langle \bd\.\rangle$.
Then $\langle \bd\.\rangle = \Delta \oplus \Lambda$ for some
free Abelian subgroup $\Lambda\in\langle \bd\.\rangle$.
Since $d_k$ and $e_k$ are indivisible in~$\Delta$,
there are $d'_{k-\?1},e'_{k-\?1}\in\Delta$ such that
$\Delta=\langle d_k,d'_{k-\?1}\rangle=\langle e_k,e'_{k-\?1}\rangle$.
After adding to  $d'_{k-\?1}$ (resp.~$e'_{k-\?1}$)
an integer multiple of $d_k$ (resp.~$e_k$),
we can assume that $0<d'_{k-\?1} < d_k$ and $0 < e'_{k-\?1} < e_k $.
Pick a basis $y_1,\dots,y_m$ of the free Abelian group~$\Lambda$.
Consider the vectors
$\bd'=(y_1,\dots,y_1,y_2,\dots,y_m,d'_{k-\?1},d_k)$,
$\be'=(y_1,\dots,y_1,y_2,\dots,y_m,e'_{k-\?1},e_k) $ in~$\RR^k$.
We have $\langle\bd\.\rangle=\langle\bd'\rangle = \langle\be'\rangle$.
After adding to the numbers $y_j$ suitable integer multiples
of $\min\.(d'_{k-\?1},e'_{k-\?1})$ and reordering them,
we can achieve that $0<y_1<\dots<y_m<d'_{k-\?1}<d_k$ and
$0<y_1<\dots<y_m<e'_{k-\?1}<e_k$, still preserving the condition 
$\langle\bd\.\rangle=\langle\bd'\rangle = \langle\be'\rangle$.

Applying Lemma~\ref{l:low-i} to the pairs $\bd,\bd'$ and $\be,\be'$, with $i=k$,
we find a low admissible path $\frak p_0$ connecting $\bd$ to~$\bd'$
and a low admissible path $\frak p_1$ connecting $\be'$ to~$\be$.
Since $y_1$ is indivisible in~$\langle \bd\.\rangle$,
we can apply  Lemma~\ref{l:low-i} to the pair $\bd',\be'$, with $i=1$,
and find a low admissible path $\frak p$ that connects $\bd'$ to~$\be'$.
Since $\bd' \le \bd$ and $\be' \le \be$, the concatenation of
$\frak p_0$, $\frak p$, and $\frak p_1$ is a low admissible path.
This completes the proof of Theorem~\ref{t:low}.
\proofend

%%%%%%%%%%%%%%%%%%%%%%%%%%%%%%%%%%%%%%%%%%%%%%%%%%%%%%%%%%%%%%%%%

\enddocument